\documentclass[small]{article}
\usepackage{amsmath,amstext,amsthm,amsfonts}
\usepackage{makeidx}
\usepackage{amssymb,latexsym}
\usepackage{txfonts, pxfonts}
\usepackage[english]{babel}
\usepackage[pdftex]{graphicx}
\numberwithin{equation}{section}
\newtheorem{theorem}{\textbf{Theorem}}[section]
\newtheorem{proposition}[theorem]{\textbf{Proposition}}
\newtheorem*{maintheorem}{\textbf{Main Theorem}}
\newtheorem{lemma}[theorem]{\textbf{Lemma}}
\newtheorem{corollary}[theorem]{Corollary}
\theoremstyle{definition}

\theoremstyle{remark}

\newenvironment{notation}[1][Notation.]{\begin{trivlist}
\item[\hskip \labelsep {\bfseries #1}]}{\end{trivlist}}

\def\e{\epsilon}

\def\R{\mathbb{R}}
\def\Rn{{\mathbb{R}}^n_+}
\def\d{\partial}

\def\esp{\Sigma}
\def\a{\alpha}
\def\b{\beta}
\def\l{\lambda}
\def\g{d}
\def\cmedia{\kappa}

\def\opinv{\mathcal{G}_{\pares}}

\def\u{u_{(\xiup,\epsilon)}}
\def\uo{u_{(0,\epsilon)}}
\def\v{v_{(\xiup,\epsilon)}}
\def\vbar{v_{(\bar{\xiup},\bar{\epsilon})}}
\def\ubar{u_{(\bar{\xiup},\bar{\epsilon})}}

\def\w{w_{(\xiup,\epsilon)}}
\def\z{z_{(\xiup,\epsilon)}}

\def\pares{(\xiup,\e)}
\def\paresbar{(\bar{\xiup},\bar{\epsilon})}
\def\esppares{\R^{n-1}\times (0,\infty)}
\def\crit {\frac{2n}{n-2}}
\def\critbordo{\frac{2(n-1)}{n-2}}
\def\conj {\frac{2n}{n+2}}
\def\conjbordo{\frac{2(n-1)}{n}}
\def\ba{\begin{align}}
\def\ea{\end{align}}
\def\bp{\begin{proof}}
\def\ep{\end{proof}}

\def\Zr{(1+x_n)^2+r^2}
\def\Zer{(\e+x_n)^2+r^2}
\def\Ze{(\e+x_n)^2+|\bar{x}-\xiup|^2}
\def\Zeo{(\e+x_n)^2+|\bar{x}|^2}

\title{Blow-up phenomena for scalar-flat metrics
on manifolds with boundary}

\author{S\'ergio de Moura Almaraz}
\date{}

\begin{document}
\maketitle

\begin{abstract}
Let $(M^n,g)$ be a compact Riemannian manifold with boundary $\d M$. This article is concerned with the set of scalar-flat metrics on $M$ which are in the conformal class of $g$ and have $\d M$ as a constant mean curvature hypersurface. We construct examples of metrics on the unit ball $B^n$, in dimensions $n\geq 25$, for which this set is noncompact. These manifolds have umbilic boundary, but they are not conformally equivalent to $B^n$.  
\end{abstract}

\section{Introduction}
Let $(M^n,g)$ be a compact Riemannian manifold with boundary $\d M$ and dimension $n\geq 3$. In 1992, J. Escobar addressed the question of finding a scalar-flat conformal metric  $\tilde{g}=u^{\frac{4}{n-2}}g$ which has $\d M$ as a constant mean curvature hypersurface. 
This problem was studied in \cite{almaraz2}, \cite{chen}, \cite{escobar3}, \cite{escobar4}, \cite{ahmedou-felli1}, \cite{coda1} and \cite{coda2}. In analytical terms, it corresponds to the existence of a positive solution to the equations
\begin{align}\label{eq:u}
\begin{cases}
\Delta_gu-c_nR_gu=0,&\text{in}\:M,
\\
\frac{\d u}{\d\eta}-d_n\cmedia_gu+Ku^{\frac{n}{n-2}}=0,&\text{on}\:\partial M,
\end{cases}
\end{align}
for some constant $K$, where $c_n=\frac{n-2}{4(n-1)}$ and $d_n=\frac{n-2}{2}$. Here, $\Delta_g$ is the Laplace-Beltrami operator, $R_g$ is the scalar curvature, $\cmedia_g$ is the mean curvature of $\d M$ and  $\eta$ is the inward unit normal vector to $\d M$. 

Escobar's question was motivated by the classical Yamabe problem, which consists of finding a conformal metric of constant scalar curvature on a given closed Riemannian manifold. This was completely solved after the works of H. Yamabe (\cite{yamabe}), N. Trudinger (\cite{trudinger}), T. Aubin (\cite{aubin1}) and R. Schoen (\cite{schoen1}). (See \cite{lee-parker} and \cite{schoen-yau3} for nice surveys on the issue.) Conformal metrics of constant scalar curvature and zero boundary mean curvature on the boundary were studied in \cite{brendle-chen}, \cite{escobar2} (see also \cite{ambrosetti-li-malchiodi} and \cite{han-li}).

The solutions to the equations (\ref{eq:u}) are the critical points of the functional 
$$
Q(u)=\frac{\int_M|d u|_g^2+c_nR_gu^2dv_g+\int_{\d M}d_n\cmedia_gu^2d\sigma_g}
{\left(\int_{\d M}|u|^\frac{2(n-1)}{n-2}d\sigma_g\right)^{\frac{n-2}{n-1}}}\,,
$$
where $dv_g$ and $d\sigma_g$ denote the volume forms of $M$ and $\d M$, respectively. 
In order to prove the existence of these solutions, Escobar introduced the conformally invariant Sobolev quotient
$$
Q(M,\d M)=\inf\{Q(u);\:u\in C^1(\bar{M}), u\nequiv 0 \:\text{on}\: \d M\}\,.
\vspace{-0.1cm}
$$

In this work we are interested in the question of whether the full set of solutions to (\ref{eq:u}) is compact. A necessary condition is that $M$ is not conformally equivalent to the standard ball $B^n$.
We point out that if the equations (\ref{eq:u}) have a solution $u>0$ with $K$ positive (resp. zero and negative), then $Q(M,\d M)$ has to be positive (resp. zero and negative). 
If $K<0$, the solution to the equations (\ref{eq:u}) is unique. If $K=0$, the equations (\ref{eq:u}) become linear and the solutions are unique up to a multiplication by a positive constant. Hence, the only interesting case is the one when $K>0$.

The problem of compactness of solutions to the equations (\ref{eq:u}) was studied by  V. Felli and M. Ould Ahmedou in the conformally flat case with umbilic boundary (\cite{ahmedou-felli1}) and in the three-dimensional case with umbilic boundary (\cite{ahmedou-felli2}). In \cite{almaraz1}, the author proved compactness for dimensions $n\geq 7$ under a generic condition. Other compactness results for similar equations were obtained by Z. Djadli, A. Malchiodi and M. Ould Ahmedou in \cite{djadli-malchiodi-ahmedou1, djadli-malchiodi-ahmedou2}, by Z. Han and Y. Li in \cite{han-li} and by M. Ould Ahmedou in \cite{ouldahmedou}.

In the case of manifolds without boundary, the question of compactness of the full set of solutions to the Yamabe equation was first raised by R. Schoen in a topics course at Stanford University in 1988.
A necessary condition is that the manifold $M^n$ is not conformally equivalent to the sphere $S^n$. 
This problem was studied in  \cite{druet1}, \cite{druet2}, \cite{li-zhang}, \cite{li-zhang2}, \cite{li-zhu2}, \cite{marques}, \cite{schoen4} and \cite{schoen-zhang} and was completely solved in a series of three papers: \cite{brendle2}, \cite{brendle-marques} and \cite{khuri-marques-schoen}. In \cite{brendle2}, S. Brendle discovered the first smooth counterexamples for dimensions $n\geq 52$ (see \cite{berti-malchiodi} for nonsmooth examples). In \cite{khuri-marques-schoen}, M. Khuri, F. Marques and R. Schoen proved compactness for dimensions $3\leq n\leq 24$. Finally, in \cite{brendle-marques}, Brendle and Marques extended the counterexamples of \cite{brendle2} to the remaining dimensions  $25\leq n\leq 51$. 

It is expected that, as in the case of manifolds without boundary, there should be a critical dimension $n_0$ such that compactness in the case of manifolds with boundary holds for $n< n_0$ and fails for $n\geq n_0$. In this work we partially answer this question by showing that compactness fails for dimensions $n\geq 25$. More precisely we prove:

\begin{maintheorem}\label{main:thm}
Let $n\geq 25$. Then there exists a smooth Riemannian metric $g$ on $B^n$ and a sequence of positive smooth functions $\{v_{\nu}\}_{\nu=1}^{\infty}$ with the following properties:

\vspace{0.1cm}
(i) $g$ is not conformally flat;

\vspace{0.1cm}
(ii) $\d B^n$ is umbilic with respect to the induced metric by $g$;

\vspace{0.1cm}
(iii) for all $\nu$, $v_{\nu}$ is a solution to the equations (\ref{eq:u}) with a constant $K>0$ and $M=B^n$; 

\vspace{0.1cm}
(iv) $Q(v_{\nu})<Q(B^n,\d B)$ for all $\nu$;

\vspace{0.1cm}
(v) $\sup_{\d B^n}v_{\nu}\to\infty$ as $\nu\to\infty$.
\end{maintheorem}

In order to prove the Main Theorem, we follow the program adopted in \cite{brendle2} and \cite{brendle-marques}. 
In Section \ref{sec:lyapunov}, we show that the problem can be reduced to finding
critical points of a certain function $\mathcal{F}_g\pares$, where $\xiup$ is a vector in $\R^{n-1}$ and
$\e$ is a positive real number. In Section \ref{sec:estim:energy}, we show that the function $\mathcal{F}_g\pares$ can be
approximated by an auxiliary function $F\pares$. In Section \ref{sec:finding}, we prove that
the function $F\pares$ has a strict local minimum point. The cases $n\geq 53$ and $25\leq n\leq 52$ are handled separately in Subsections \ref{subsec:case53} and \ref{subsec:case25} respectively.
Finally, in Section \ref{sec:proof}, we use a perturbation argument to construct critical
points of the function $\mathcal{F}_g\pares$ and prove the non-compactness theorem.

\begin{notation}
Throughout this work we will make use of the index notation for tensors. We will adopt the summation convention whenever confusion is not possible and use indices $1\leq i,i,j,k,l,m,p,q,r,s\leq n-1$ and $1\leq a,b,c,d\leq n$. 
We also define constants $c_n=\frac{n-2}{4(n-1)}$ and $d_n=\frac{n-2}{2}$.

We will denote by $\Delta_g$ the Laplace-Beltrami operator. 
The volume forms of $M$ and $\d M$ will be denoted by $dv_g$ and $d\sigma_g$, respectively. By $\eta$ we will denote the inward unit normal vector to $\d M$.
The scalar curvature will be denoted by $R_g$ , the second fundamental form of $\d M$ by $\pi_{kl}$ and the mean curvature, $\frac{1}{n-1}tr (\pi_{kl})$, by $\cmedia_g$.

By $\Rn$ we will denote the half-space $\{x=(x_1,...,x_n)\in \R^n;\:x_n\geq 0\}$. If $x\in\Rn$ we set $\bar{x}=(x_1,...,x_{n-1},0)\in\d\Rn\cong \R^{n-1}$. For any $x_0\in\Rn$ we set $B^+_r(x_0)=\{x\in\Rn\,;\:|x-x_0|<r\}$.
The n-dimensional sphere of radius $r$ in $\R^{n+1}$ will be denoted by $S_r^n$ and 
$\sigma_{n}$ will denote the area of the n-dimensional unit sphere $S^n_1$.   
\end{notation}
{\bf{Acknowledgements.}}
I would like to thank Prof. Fernando C. Marques for his comments and interest in this work.

\section{Lyapunov-Schmidt reduction}\label{sec:lyapunov}
Given a pair $(\xiup,\e)\in \esppares$ we set 
$$
\u(x)=\left(\frac{\e}{\Ze}\right)^{\frac{n-2}{2}}\,,\:\:\:\:\text{for}\:x\in\Rn\,.
$$
Observe that $\u$ satisfies
\begin{align}\label{eq:U}
\begin{cases}
\Delta\u=0,&\text{in}\:\Rn,
\\
\frac{\d}{\d x_n}\u+(n-2)\u^{\frac{n}{n-2}}=0,&\text{on}\:\partial \Rn,
\end{cases}
\end{align}
and 
\begin{equation}\label{eq:u:Q}
\int_{\d\Rn}\u^{\frac{2(n-1)}{n-2}}=\left(\frac{Q(B^n,\d B)}{n-2}\right)^{n-1}\,.
\end{equation}

Let us define 
$$
\phi_{(\xiup,\e,n)}(x)=\left(\frac{\e}{\Ze}\right)^{\frac{n}{2}}\frac{\e^2-x_n^2-|\bar{x}-\xiup|^2}{\Ze}
$$
and
$$
\phi_{(\xiup,\e,k)}(x)=\left(\frac{\e}{\Ze}\right)^{\frac{n}{2}}\frac{2\e(x_k-\xiup_k)}{\Ze}
$$
for $x\in\Rn$ and $k=1,...,n-1$. Observe that 
$$
\phi_{(\xiup,\e,n)}(x)\cdot(\Ze)=-\frac{2\e^{2}}{n-2}\frac{\d}{\d\e}\u(x)\,,
$$
$$
\phi_{(\xiup,\e,k)}(x)\cdot(\Ze)=\frac{2\e^{2}}{n-2}\frac{\d}{\d\xiup_k}\u(x)\,,
$$
for $k=1,...,n-1$, and that $\|\phi_{(\xiup,\e,a)}\|_{L^{\conjbordo}(\d\Rn)}$ is independent of $(\e,\xiup)\in\esppares$, for any $a=1,...,n$.  

We also set
$$
\esp=\left\{w\in L^{\crit}(\Rn)\cap L^{\critbordo}(\d\Rn)\cap H^1_{loc}(\Rn)\,;\:\int_{\Rn}|dw|^2<\infty\right\}\,,
$$
$$
\esp_{\pares}=\left\{w\in\esp\,;\:\int_{\d\Rn}\phi_{(\xiup,\e,a)}w=0\,,\:a=1,...,n \right\}
$$
and  $\|w\|_{\esp}=\left(\int_{\Rn}|dw|^2\right)^{\frac{1}{2}}$ for $w\in\esp$. Observe that $\u\in\esp_{\pares}$ for each $\pares\in\esppares$. By Sobolev's inequality, there exists $K=K(n)>0$ such that
\begin{equation}\label{estim:K}
\left(\int_{\Rn}|w|^{\frac{2n}{n-2}}\right)^{\frac{n-2}{n}}
+\left(\int_{\d\Rn}|w|^{\frac{2(n-1)}{n-2}}\right)^{\frac{n-2}{n-1}}
\leq K\int_{\Rn}|dw|^2
\end{equation}
for all $w\in\esp$.

In what follows in this section we are going to find, for each pair $\pares\in\esppares$, a function $\v\in\esp$ which is an approximate weak solution to a Yamabe-type problem (\ref{eq:u}) on $\Rn$. Then we will show that $\v$ is in fact a classical solution to this problem whenever $\pares$ is a critical point of a certain energy function defined on $\esppares$.  

The following result is Proposition 26 of \cite{brendle2} and will be used throughout this work:
\begin{lemma}\label{curv_esc}
Suppose that we express the Riemannian metric $g$ as $g=exp(h)$, where $h$ is a trace-free symmetric  two-tensor defined on $\Rn$ and satisfying $|h(x)|\leq 1$ for any $x\in\Rn$. Then there exists $C=C(n)>0$ such that
$$
\left|R_g-
\left\{
\d_a\d_bh_{ab}-\d_a(h_{ac}\d_bh_{bc})
+\frac{1}{2}\d_ah_{ac}\d_bh_{bc}
-\frac{1}{4}\d_ch_{ab}\d_ch_{ab}
\right\}
\right|
\leq
C|h|^2|\d^2h|+C|h||\d h|^2\,.
$$
\end{lemma}
\begin{notation}
In this section we suppose that $g$ is a Riemannian metric on $\Rn$ expressed as $g=exp(h)$, where $h$ is a trace-free symmetric two-tensor satisfying $h(x)=0$ for any $|x|\geq 1$.  
\end{notation}
\begin{proposition}\label{Propo1}
If $|h(x)|+|\d h(x)|+|\d^2h(x)|\leq \a\leq 1$ for any $x\in\Rn$, then there exists $C=C(n)>0$ such that
$$
\left\|\Delta_g\u-c_nR_g\u\right\|_{L^{\frac{2n}{n+2}}(\Rn)}
+\left\|d_n\cmedia_g\u\right\|_{L^{\frac{2(n-1)}{n}}(\d\Rn)}
\leq C\a
$$
for all pairs $(\xiup,\e)\in \esppares$.
\end{proposition}
\bp
It follows from the pointwise estimates
$$
\left|\Delta_g\u-c_nR_g\u\right|\leq C\left\{|h||\d^2\u|+|\d h||\d\u|+(|\d^2h|+|\d h|^2)|\u|\right\}
$$
and $|d_n\cmedia_g\u|
\leq C|\d h||\u|$
that
\ba
&\left\|\Delta_g\u-c_nR_g\u\right\|_{L^{\frac{2n}{n+2}}(\Rn)}\notag
\\
&\hspace{1cm}\leq C\left\{\|h\|_{L^{\infty}(\Rn)}\|\d^2\u\|_{L^{\frac{2n}{n+2}}(\Rn)}+\|\d h\|_{L^{n}(\Rn)}\|\d\u\|_{L^{2}(\Rn)}\right\}\notag
\\
&\hspace{1.5cm}+C\left\{\|\d^2h\|_{L^{\frac{n}{2}}(\Rn)}+\|\d h\|^2_{L^{n}(\Rn)}\right\}\|u\|_{L^{\crit}(\Rn)}\notag
\end{align}
and
$$
\|d_n\cmedia_g\u\|_{L^{\frac{2(n-1)}{n}}(\d\Rn)}
\leq C\|\d h\|_{L^{n-1}(\d\Rn)}\|\u\|_{L^{\frac{2(n-1)}{n-2}}(\d\Rn)}\,.\notag
$$
From this the result follows.
\ep
\begin{lemma}\label{lemma:eq:conf}
Let $B^n=B^n_{1/2}(0,...,0,-\frac{1}{2})\subset\R^n$ be the ball with radius $\frac{1}{2}$ and center $(0,...,0,-\frac{1}{2})$. Let $z_1,...,z_n$ be the coordinates of $B^n$ taken with center $(0,...,0,-\frac{1}{2})$. For each pair $\pares\in\esppares$, the expression
$$
\mathcal{C}_{\pares}(x)=\frac{\e\,(x_1-\xiup_1,...,x_{n-1}-\xiup_{n-1},x_n+\e)}{|\bar{x}-\xiup|^2+(x_n+\e)^2}
+(0,...,0,-1)
$$
defines a conformal equivalence 
$$
\mathcal{C}_{\pares}:\Rn\to B^n\backslash\{(0,...,0,-1)\}
$$ 
that satisfies $\mathcal{C}_{\pares}^{\:\:*}\delta_{B^n}=\u^{\frac{4}{n-2}}\delta$, where $\delta_{B^n}$ is the Euclidean metric on $B^n$ and $\delta$ is the Euclidean metric on $\Rn$. For any smooth function $f$ on $\Rn$, we have 
\begin{equation}\label{conf:L:bola}
\Delta_{B^n}\tilde{u}_{\pares}=\u^{-\frac{n+2}{n-2}}\Delta f
\end{equation}
and
\begin{equation}\label{conf:B:bola}
\frac{\d}{\d\eta}\tilde{u}_{\pares}-(n-2)\tilde{u}_{\pares}=\u^{-\frac{n}{n-2}}\frac{\d f}{\d x_n}\,,
\end{equation} 
where $\tilde{u}_{\pares}=(f\u^{-1})\circ\mathcal{C}_{\pares}^{\:-1}$. Moreover,
\begin{equation}\label{expr:zn}
z_n\circ\mathcal{C}_{\pares}=-\frac{\e}{n-2}\u^{-1}\frac{\d}{\d\e}\u
=\frac{1}{2}\u^{-\frac{n}{n-2}}\phi_{(\xiup,\e,n)}
\end{equation}
and
\begin{equation}\label{expr:zk}
z_k\circ\mathcal{C}_{\pares}=\frac{\e}{n-2}\u^{-1}\frac{\d}{\d\xiup_k}\u
=\frac{1}{2}\u^{-\frac{n}{n-2}}\phi_{(\xiup,\e,k)}\,,
\:\:\:\:k=1,...,n-1\,.
\end{equation}
\end{lemma}
\bp
These are direct computations. The assertions (\ref{conf:L:bola}) and (\ref{conf:B:bola}) follow from the following properties of the conformal operators $L_g=\Delta_g-c_nR_g$ and $B_g=\frac{\d}{\d\eta}-d_n\cmedia_g$:
\begin{equation}\label{conf:proper}
L_{u^{\frac{4}{n-2}}g}(fu^{-1})=u^{-\frac{n+2}{n-2}}L_g f 
\:\:\:\:\:\text{and}\:\:\:\:\:
B_{u^{\frac{4}{n-2}}g}(fu^{-1})=u^{-\frac{n}{n-2}}B_g f\,. 
\end{equation} 
\ep
\begin{lemma}\label{lema:bola}
There exists $\theta=\theta(n)>0$ such that
$$
\int_{B^n}|dw|^2-2\int_{\d B^n}w^2-2\theta\left(\int_{B^n}|dw|^2+(n-2)\int_{\d B^n}w^2\right)+\frac{4}{\theta}\left(\int_{\d B^n}w\right)^2\geq 0
$$
for any $w\in H^1(B^n)$ such that $w\perp_{L^2(\d B^n)} \{z_1,...,z_{n}\}$. Here, we are following the notations of Lemma \ref{lemma:eq:conf}.
\end{lemma}
\bp
First we fix $0\nequiv w\in H^1(B^n)$ such that $w\perp_{L^2(\d B^n)} \{1,z_1,...,z_{n}\}$. Since 
$$
\inf \left\{ 
\frac{\int_{B^n}|d\psi|^2}{\int_{\d B^n}\psi^2}\,,\:\text{such that}\: 
\psi\in H^1(B^n),\,\psi\nequiv 0\:{on}\:\d B^n\:\text{and}\: \psi\perp_{L^2(\d B^n)} 1 
\right\}=2
$$
and this infimum is realized only by the functions $z_1,...,z_n$, we see that 
$$
\int_{B^n}|dw|^2-2\int_{\d B^n}w^2>0\,.
$$ 
Hence,
\begin{equation}\label{lema:bola:1}
\int_{B^n}|dw|^2-2\int_{\d B^n}w^2\geq 2\theta\left(\int_{B^n}|dw|^2+(n-2)\int_{\d B^n}w^2\right)
\end{equation}
holds for any $\theta>0$ satisfying
$$
\theta\leq \theta(w)=\frac{1}{2}\frac{\int_{B^n}|dw|^2-2\int_{\d B^n}w^2}{\int_{B^n}|dw|^2+(n-2)\int_{\d B^n}w^2}
$$ 
and the equality is realized by $\theta=\theta(w)$.

We claim that there exists $\theta_0>0$ such that $\theta(w)\geq\theta_0$ for any $w\in H^1(B^n)$ satisfying $w\perp_{L^2(\d B^n)} \{1,z_1,...,z_{n}\}$. Suppose by contradiction this is not true. Thus we can choose a sequence $\{w_j\}_{j=1}^{\infty}\subset H^1(B^n)$ such that $w_j\perp_{L^2(\d B^n)}\{1,z_1,...,z_{n}\}$ and $\theta(w_j)\to 0$ as $j\to \infty$. Hence 
\begin{equation}\notag
\int_{B^n}|dw_j|^2-2\int_{\d B^n}w_j^2
= 2\theta(w_j)\left(\int_{B^n}|dw_j|^2+(n-2)\int_{\d B^n}w_j^2\right)
\end{equation}
holds and we can assume that $\int_{B^n}|dw_j|^2=1$ for any $j$. Thus, $\int_{\d B^n}w_j^2\leq\frac{1}{2}$ for all $j$ and we can suppose that $w_j\rightharpoonup w_0$ in $H^1(B^n)$ for some $w_0$. Since $H^1(B^n)$ is compactly imbedded in $L^{2}(\d B^n)$, we know that $w_0\perp_{L^2(\d B^n)} \{1,z_1,...,z_{n}\}$. Let us first assume that $w_0\nequiv 0$.
We set
$$
\b=\int_{B^n}|dw_0|^2-2\int_{\d B^n}w_0^2>0\,.
$$
Since $\liminf_{i\to\infty}\int_{B^n}|dw_j|^2\geq \int_{B^n}|dw_0|^2$ and $\lim_{i\to\infty}\int_{\d B^n}w_j^2= \int_{\d B^n}w^2$, we can assume that
$$
\int_{B^n}|dw_j|^2-2\int_{\d B^n}w_j^2\geq\frac{\b}{2}\:\:\:\:\text{for all}\:j\,.
$$ 
On the other hand,
$$
\frac{\b}{n}
\left\{\int_{B^n}|dw_j|^2+(n-2)\int_{\d B^n}w_j^2\right\}
\leq \frac{\b}{2}\,,
$$
since $\int_{B^n}|dw_j|^2=1$ and $\int_{\d B^n}w_j^2\leq\frac{1}{2}$. 
Hence,
\ba
2\theta(w_j)\left(\int_{B^n}|dw_j|^2+(n-2)\int_{\d B^n}w_j^2\right)
&=\int_{B^n}|dw_j|^2-2\int_{\d B^n}w_j^2\notag
\\
&\geq\frac{\b}{n}\left(\int_{B^n}|dw_j|^2+(n-2)\int_{\d B^n}w_j^2\right)\,.\notag
\end{align}
which implies that $2\theta(w_j)\geq \frac{\b}{n}$ for all $j$ and contradicts the fact that $\theta(w_j)\to 0$. 

Thus we must have $w_0\equiv 0$, which implies that $\int_{\d B^n}w_j^2\to 0$ as $j\to\infty$. Then, if we set $\tilde{w}_j=\left(\int_{\d B^n}w_j^2\right)^{-\frac{1}{2}}w_j$, we have
$\tilde{w}_j\rightharpoonup\tilde{w}_0$ in $H^1(B^n)$, for some $\tilde{w}_0$. Moreover, 
$$
0=\lim_{j\to\infty}\int_{B^n}|d\tilde{w}_j|^2\geq \int_{B^n}|d\tilde{w}_0|^2
$$
and
$$
\int_{\d B^n}\tilde{w}_j^2=1=\int_{\d B^n}\tilde{w}_0^2\,.
$$
From this we conclude that $\tilde{w}_0\equiv \text{const}\neq 0$, which contradicts the fact that $\tilde{w}_0\perp_{L^2(\d B^n)} 1$. This proves that there exists $\theta_0>0$ such that $\theta(w)\geq\theta_0$ for any $w\in H^1(B^n)$ satisfying $w\perp_{L^2(\d B^n)} \{1,z_1,...,z_{n}\}$. In particular, (\ref{lema:bola:1}) holds, with $\theta=\theta_0$, for any $w\in H^1(B^n)$ satisfying $w\perp_{L^2(\d B^n)} \{1,z_1,...,z_{n}\}$.

Now, let $w\in H^1(B^n)$ satisfy $w\perp_{L^2(\d B^n)} \{z_1,...,z_{n}\}$. We write $w=w_1+b$ where $b$ is a constant and $w_1\perp_{L^2(\d B^n)} 1$. Then we have
\ba
&\int_{B^n}|dw|^2-2\int_{\d B^n}w^2-2\theta_0\left(\int_{B^n}|dw|^2+(n-2)\int_{\d B^n}w^2\right)+\frac{4}{\theta_0}\left(\int_{\d B^n}w\right)^2\notag
\\
&\hspace{1cm}=\int_{B^n}|dw_1|^2-2\int_{\d B^n}w_1^2
-2\theta_0\left(\int_{B^n}|dw_1|^2+(n-2)\int_{\d B^n}w_1^2\right)\notag
\\
&\hspace{2cm}-2(1+(n-2)\theta_0)\int_{\d B^n}b^2+\frac{4}{\theta_0}\left(\int_{\d B^n}b\right)^2\notag
\\
&\hspace{1cm}\geq \left(\frac{4}{\theta_0}-2-2(n-2)\theta_0\right)\int_{\d B^n}b^2\notag
\end{align}
Choosing $\theta_0$ smaller if necessary, we can suppose that $\frac{4}{\theta_0}-2-2(n-2)\theta_0>0$ and the result follows.
\ep 
\begin{proposition}\label{Propo2}
There exists $\theta=\theta(n)>0$ such that
$$
\int_{\Rn}|dw|^2-n\int_{\d\Rn}\u^{\frac{2}{n-2}}w^2\geq 2\theta|w|^2_{\esp}-\frac{4}{\theta}\left(\int_{\d\Rn}\u^{\frac{n}{n-2}}w\right)^2
$$
for all $w\in\esp_{(\xiup,\e)}$ and any pair $\pares\in\esppares$.
\end{proposition}
\bp
Let $w\in\esp_{\pares}$ and set $\bar{w}=(w\u^{-1})\circ \mathcal{C}_{\pares}^{\:-1}$. 
Using the fact that $C_0^{\infty}(\Rn)$ is dense in $\esp$ with respect to the norms $\|\cdot\|_{\esp}$, $\|\cdot\|_{L^{\crit}(\Rn)}$ and $\|~\cdot~\|_{L^{\critbordo}(\d\Rn)}$, it is easy to see that we can assume that $\bar{w}\in H^1(B^n)$.   
It follows from the expressions (\ref{expr:zn}) and (\ref{expr:zk}) that
$$
\int_{\d B^n}\bar{w}\,z_n
=\frac{1}{2}\int_{\d\Rn}w\,\phi_{(\xiup,\e,n)}=0
$$
and
$$
\int_{\d B^n}\bar{w}\,z_k
=\frac{1}{2}\int_{\d\Rn}w\,\phi_{(\xiup,\e,k)}=0\,,
\:\:\:\:k=1,...,n-1\,.
$$
Then, according to Lemma \ref{lema:bola}, we have
\begin{equation}\label{Propo2:1}
\int_{B^n}|d\bar{w}|^2-2\int_{\d B^n}\bar{w}^2-2\theta\left(\int_{B^n}|d\bar{w}|^2+(n-2)\int_{\d B^n}\bar{w}^2\right)+\frac{4}{\theta}\left(\int_{\d B^n}\bar{w}\right)^2\geq 0\,.
\end{equation}
Hence, using the formulas (\ref{conf:L:bola}) and (\ref{conf:B:bola}), we ealisy see that    
$$
\int_{B^n}|d\bar{w}|^2-2\int_{\d B^n}\bar{w}^2
=\int_{\Rn}|dw|^2-n\int_{\d\Rn}\u^{\frac{2}{n-2}}w^2\,,
$$
$$
\int_{B^n}|d\bar{w}|^2+(n-2)\int_{\d B^n}\bar{w}^2
=\int_{\Rn}|dw|^2
\:\:\:\:\:\text{and}\:\:\:\:\:
\int_{\d B^n}\bar{w}=\int_{\d\Rn}\u^{\frac{n}{n-2}}w\,.
$$
Now the result follows from substituting these last three equations in (\ref{Propo2:1}).
\ep 
\begin{corollary}\label{Corol3}
Let $K$ be as in (\ref{estim:K}) and $\theta$ be as in Proposition \ref{Propo2}. Then there exists $0<\a_0=\a_0(n)\leq 1$ such that, whenever $|h(x)|+|\d h(x)|+|\d^2h(x)|\leq\a_0$ for all $x\in\Rn$, we have
\begin{equation}\label{Corol3:4}
\left(\int_{\Rn}|w|^{\frac{2n}{n-2}}\right)^{\frac{n-2}{n}}
+\left(\int_{\d\Rn}|w|^{\frac{2(n-1)}{n-2}}\right)^{\frac{n-2}{n-1}}
\leq
2K\int_{\Rn}\left(|dw|_g^2+c_nR_gw^2\right)+2K\int_{\d\Rn}d_n\cmedia_g w^2
\end{equation}
for all $w\in\esp$ and
\begin{equation}\label{Corol3:3}
\int_{\Rn}(|dw|_g^2+c_nR_gw^2)
+\int_{\d\Rn}\left(d_n\cmedia_g w^2-n\u^{\frac{2}{n-2}}w^2
\right)
\geq\frac{\theta}{2}\|w\|_{\esp}^2-\frac{1}{\theta}A(w)^2
\end{equation}
for all $w\in\esp_{(\xiup,\e)}$ and any pair $\pares\in\esppares$. Here,
$$
A(w)=\int_{\Rn}\left(\Delta_g\u-c_nR_g\u\right)w
+\int_{\d\Rn}\left(-d_n\cmedia_g\u +2u^{\frac{n}{n-2}}\right)w\,.
$$
\end{corollary}
\bp
Let us first prove the estimate (\ref{Corol3:3}). Observe that
$$
\int_{\Rn}\left(\Delta_g\u-c_nR_g\u\right)w
\geq -\left\|\Delta_g\u-c_nR_g\u\right\|_{L^{\conj}(\Rn)}
\left\|w\right\|_{L^{\crit}(\Rn)}
$$
and
\ba
&\int_{\d\Rn}\left(-d_n\cmedia_g\u +2u^{\frac{n}{n-2}}\right)w\notag
\\
&\hspace{2cm}\geq -\left\|d_n\cmedia_g\u\right\|_{L^{\conjbordo}(\d\Rn)}
\left\|w\right\|_{L^{\critbordo}(\d\Rn)}
+2\int_{\d\Rn}\u^{\frac{n}{n-2}}w\notag
\end{align}
Hence, by Proposition \ref{Propo1} and inequality (\ref{estim:K}) we have
$$
A(w)\geq -C\a_0\|w\|_{\esp}+2\int_{\d\Rn}\u^{\frac{n}{n-2}}w\,.
$$
Choosing $\a_0$ small this implies
$$
A(w)^2\geq 4\left(\int_{\d\Rn}\u^{\frac{n}{n-2}}w\right)^2-\theta^2\|w\|^2_{\esp}
$$
which, together with Proposition \ref{Propo2}, gives
\begin{equation}\label{Corol3:1}
\int_{\Rn}|dw|^2-\int_{\d\Rn}n\u^{\frac{2}{n-2}}w^2\geq\theta\|w\|^2_{\esp}-\frac{1}{\theta}A(w)^2\,.
\end{equation}
On the other hand,
\ba
&\int_{\Rn}(|dw|_g^2+c_nR_gw^2)+\int_{\d\Rn}(d_n\cmedia_g w^2-n\u^{\frac{2}{n-2}}w^2)\notag
\\
&=\int_{\Rn}|dw|^2-\int_{\d\Rn}n\u^{\frac{2}{n-2}}w^2
+\int_{\Rn}\left\{(g^{ij}-\delta^{ij})\d_iw\d_jw+c_nR_gw^2\right\}+\int_{\d\Rn}d_n\cmedia_g w^2\,.\notag
\end{align}
The fact that $h(x)=0$ for $|x|\geq 1$ and (\ref{estim:K}) imply that
\ba\label{Corol3:5}
&\int_{\Rn}\left\{(g^{ij}-\delta^{ij})\d_iw\d_jw+c_nR_gw^2\right\}+\int_{\d\Rn}d_n\cmedia_g w^2
\\
&\hspace{2cm}\leq C\a_0\|w\|^2_{\esp}+C\a_0\|w\|_{L^{\crit}(\Rn)}^2+C\a_0\|w\|_{L^{\critbordo}(\d\Rn)}^2\notag
\\
&\hspace{2cm}\leq C\a_0(1+K)\|w\|^2_{\esp}\,.\notag
\end{align}
Hence,
\ba\label{Corol3:2}
&\int_{\Rn}(|dw|_g^2+c_nR_gw^2)+\int_{\d\Rn}(d_n\cmedia_g w^2-n\u^{\frac{2}{n-2}}w^2)
\\
&\hspace{2cm}\geq \int_{\Rn}|dw|^2-\int_{\d\Rn}n\u^{\frac{2}{n-2}}w^2-C\a_0(1+K)\|w\|^2_{\esp}\,.\notag
\end{align}
Now the result follows from the inequalities (\ref{Corol3:1}) and (\ref{Corol3:2}), choosing $\a_0$ small.
The estimate (\ref{Corol3:4}) follows easily from the inequalities (\ref{estim:K}) and (\ref{Corol3:5}). 
\ep
\begin{proposition}\label{Propo4}
Suppose that $|h(x)|+|\d h(x)+|\d^2 h(x)|\leq \a_0$ for all $x\in\Rn$, where $\a_0$ is the constant obtained in Corollary \ref{Corol3}.
Given any pair $(\xiup,\e)\in \esppares$ and functions $f\in L^{\frac{2n}{n+2}}(\Rn)$ and $\bar{f}\in L^{\frac{2(n-1)}{n}}(\d\Rn)$ there exists a unique $w\in\esp_{(\xiup,\e)}$ such that
\begin{equation}\label{Propo4:1}
\int_{\Rn}\left(<dw,d\psi>_g+c_nR_gw\psi\right)
+\int_{\d\Rn}\left(d_n\cmedia_g w\psi-n\u^{\frac{2}{n-2}}w\psi\right)
=\int_{\Rn}f\psi+\int_{\d\Rn}\bar{f}\psi
\end{equation}
for all $\psi\in\esp_{(\xiup,\e)}$. Moreover, there exists $C=C(n)>0$ such that
$$
\|w\|_{\esp}\leq C\|f\|_{L^{\frac{2n}{n+2}}(\Rn)}+C\|\bar{f}\|_{L^{\frac{2(n-1)}{n}}(\d\Rn)}\,.
$$
\end{proposition}
\bp
Let us first prove the existence part. Following the notations of Corollary \ref{Corol3}, we define the funcional
$$
T(w)=\int_{\Rn}(|dw|_g^2+c_nR_gw^2-2fw)
+\int_{\d\Rn}(d_n\cmedia_g w^2-n\u^{\frac{2}{n-2}}w^2-2\bar{f}w)
+\frac{1}{\theta}A(w)^2
$$
for $w\in\esp_{\pares}$. Hence
\ba
dT_w(\psi)
&=2\int_{\Rn}\left(<dw,d\psi>_g+c_nR_gw\psi-f\psi\right)
+2\int_{\d\Rn}\left(d_n\cmedia_g w\psi-n\u^{\frac{2}{n-2}}w\psi-\bar{f}\psi\right)\notag
\\
&\hspace{0.5cm}+\frac{2}{\theta}A(w)A(\psi)\,.\notag
\end{align}

It follows from the identity (\ref{Corol3:3}) that
\ba
T(w)&\geq \frac{\theta}{2}\|w\|^2_{\esp}-2\int_{\Rn}fw-2\int_{\d\Rn}\bar{f}w\notag
\\
&\geq \frac{\theta}{2}\|w\|^2_{\esp}-2\|f\|_{L^{\conj}(\Rn)}\|w\|_{L^{\crit}(\Rn)}
-2\|\bar{f}\|_{L^{\conjbordo}(\d\Rn)}\|w\|_{L^{\critbordo}(\d\Rn)}\notag
\\
&\geq \frac{\theta}{4}\|w\|^2_{\esp}-C\left(\|f\|_{L^{\conj}(\Rn)}^2+\|\bar{f}\|_{L^{\conjbordo}(\d\Rn)}^2\right)\notag
\end{align} 
where in the last inequality we used the estimate (\ref{estim:K}). So, $T$ is bounded below and by a standard argument we can find a minimizer $w_0$ for $T$ over all functions in $\esp_{\pares}$.
Now, integrating by parts we see that
$$
\int_{\Rn}\left(<d\u,d\psi>_g+c_nR_g\u\psi\right)
+\int_{\d\Rn}\left(d_n\cmedia_g\u\psi-n\u^{\frac{n}{n-2}}\psi\right)
=-A(\psi)\,,
$$
holds for all $\psi\in C_0^{\infty}(\Rn)$. Since this space is dense in $\esp$ with respect to the norms $\|\cdot\|_{\esp}$, $\|\cdot\|_{L^{\crit}(\Rn)}$ and $\|~\cdot~\|_{L^{\critbordo}(\d\Rn)}$, this identity holds for all $\psi\in\esp$.
Hence, the function $w=w_0-\frac{1}{\theta}A(w_0)\,\u$ satisfies (\ref{Propo4:1}) for all $\psi\in\esp_{\pares}$, proving the existence part.

In order to prove the uniqueness part, suppose that $w\in\esp_{(\xiup,\e)}$ satisfies (\ref{Propo4:1}) for all $\psi\in\esp_{(\xiup,\e)}$. In particular,
$$
\int_{\Rn}\left(|dw|_g^2+c_nR_gw^2\right)
+\int_{\d\Rn}\left(d_n\cmedia_g w^2-n\u^{\frac{2}{n-2}}w^2\right)
=\int_{\Rn}fw+\int_{\d\Rn}\bar{f}w
$$
and
\begin{align}
-A(w)
&=\int_{\Rn}\left(<dw,d\u>_g+c_nR_gw\u\right)
+\int_{\d\Rn}\left(d_n\cmedia_g w\u-n\u^{\frac{n}{n-2}}w\right)\notag
\\
&=\int_{\Rn}f\u+\int_{\d\Rn}\bar{f}\u\,,\notag
\end{align}
since $\u\in\esp_{\pares}$. Then (\ref{Corol3:3}) implies
\ba
\frac{\theta}{2}\|w\|_{\esp}^2
&\leq \int_{\Rn}(|dw|_g^2+c_nR_gw^2)
+\int_{\d\Rn}\left(d_n\cmedia_g w^2-n\u^{\frac{2}{n-2}}w^2
\right)+\frac{1}{\theta}A(w)^2\notag
\\
&=\int_{\Rn}fw+\int_{\d\Rn}\bar{f}w
+\frac{1}{\theta}\left(\int_{\Rn}f\u+\int_{\d\Rn}\bar{f}\u\right)^2\notag
\\
&\leq \left\{\|w\|_{L^{\crit}(\Rn)}
+\frac{2}{\theta}\|\u\|^2_{L^{\crit}(\Rn)}\|f\|_{L^{\conj}(\Rn)}\right\}
\|f\|_{L^{\conj}(\Rn)}\notag
\\
&\hspace{1cm}+\left\{\|w\|_{L^{\critbordo}(\d\Rn)}
+\frac{2}{\theta}\|\u\|^2_{L^{\critbordo}(\d\Rn)}\|\bar{f}\|_{L^{\conjbordo}(\d\Rn)}\right\}
\|\bar{f}\|_{L^{\conjbordo}(\d\Rn)}\notag
\\
&\leq \left\{K^{\frac{1}{2}}\|w\|_{\esp}+\frac{2}{\theta}\|\u\|^2_{L^{\crit}(\Rn)}\|f\|_{L^{\conj}(\Rn)}\right\}
\|f\|_{L^{\conj}(\Rn)}\notag
\\
&\hspace{1cm}+\left\{K^{\frac{1}{2}}\|w\|_{\esp}
+\frac{2}{\theta}\|\u\|^2_{L^{\critbordo}(\d\Rn)}\|\bar{f}\|_{L^{\conjbordo}(\d\Rn)}\right\}
\|\bar{f}\|_{L^{\conjbordo}(\d\Rn)}\,.\notag
\end{align}
Hence,
\ba
\frac{\theta}{4}\|w\|_{\esp}^2
&\leq \left\{\frac{K}{\theta}+\frac{2}{\theta}\|\u\|^2_{L^{\crit}(\Rn)}\right\}\|f\|^2_{L^{\conj}(\Rn)}\notag
\\
&\hspace{1cm}+\left\{\frac{K}{\theta}+\frac{2}{\theta}\|\u\|^2_{L^{\critbordo}(\d\Rn)}\right\}
\|\bar{f}\|^2_{L^{\conjbordo}(\d\Rn)}\notag
\end{align}
and the result follows.
\ep
\begin{proposition}\label{Propo5}
Let $\a_0$ be the constant obtained in Corollary \ref{Corol3}.
There is a constant $\a_1=\a_1(n)$, $0<\a_1\leq\a_0$, with the following property:
if $|h(x)|+|\d h(x)|+|\d^2h(x)|\leq \a_1$ for all $x\in\Rn$, given any pair $(\xiup,\e)\in\esppares$ there exists a unique $\v\in\esp$ such that $\v-\u\in\esp_{(\xiup,\e)}$ and
\begin{equation}\notag
\int_{\Rn}\left(<d\v,d\psi>_g+c_nR_g\v\psi\right)
+\int_{\d\Rn}\left(d_n\cmedia_g\v\psi-(n-2)|\v|^{\frac{2}{n-2}}\v\psi\right)=0
\end{equation}
for all $\psi\in\esp_{(\xiup,\e)}$. Moreover, there exists $C=C(n)>0$ such that
\begin{equation}\label{Propo5:1}
\|\v-\u\|_{\esp}\leq C\|\Delta_g\u-c_nR_g\u\|_{L^{\frac{2n}{n+2}}(\Rn)}+C\|d_n\cmedia_g\u\|_{L^{\frac{2(n-1)}{n}}(\d\Rn)}
\end{equation}
In particular, $\v\nequiv 0$.
\end{proposition}
\bp
Using Proposition \ref{Propo4} we can define
$$
\opinv:L^{\conj}(\Rn)\times L^{\conjbordo}(\d\Rn)\longrightarrow \esp_{\pares}
$$
by $\opinv(f,\bar{f})=w$, where $w\in\esp_{\pares}$ satisfies (\ref{Propo4:1}) for all $\psi\in\esp_{\pares}$. Hence, there exists $C=C(n)$ such that
\begin{equation}\label{Propo5:2}
\|\opinv(f,\bar{f})\|_{\esp}\leq C\|f\|_{L^{\conj}(\Rn)}+C\|\bar{f}\|_{L^{\conjbordo}(\d\Rn)}\,.
\end{equation}
We define a nonlinear mapping $\Phi_{\pares}(w):\esp_{\pares}\to \esp_{\pares}$ by
$$
\Phi_{\pares}(w)=\opinv(f_{\pares},\bar{f}_{(\xiup,\e,w)})
$$
where
$$
f_{\pares}=\Delta_g\u-c_nR_g\u
$$
and
$$
\bar{f}_{(\xiup,\e,w)}=-d_n\cmedia_g\u+(n-2)\left\{|\u+w|^{\frac{2}{n-2}}(\u+w)
-\u^{\frac{n}{n-2}}-\frac{n}{n-2}\u^{\frac{2}{n-2}}w\right\}\,.\notag
$$
It follows from Proposition \ref{Propo1} and the inequality (\ref{Propo5:2}) that
$\|\Phi_{\pares}(0)\|_{\esp}\leq C\a_1$.
Since
\ba
&\left| |\u+w|^{\frac{2}{n-2}}(\u+w)-|\u+\tilde{w}|^{\frac{2}{n-2}}(\u+\tilde{w})
-\frac{n}{n-2}\u^{\frac{2}{n-2}}(w-\tilde{w})\right|\notag
\\
&\hspace{1cm}\leq C\left(|w|^{\frac{2}{n-2}}+|\tilde{w}|^{\frac{2}{n-2}}\right)|w-\tilde{w}|\,,\notag
\end{align}
we have
\ba
&\|\Phi_{\pares}(w)-\Phi_{\pares}(\tilde{w})\|_{\esp}\notag
\\
&\hspace{1cm}\leq C\left\|\left(|w|^{\frac{2}{n-2}}
+|\tilde{w}|^{\frac{2}{n-2}}\right)(w-\tilde{w})\right\|_{L^{\conjbordo}(\d\Rn)}\notag
\\
&\hspace{1cm}\leq C\left\{\|w\|_{L^{\critbordo}(\d\Rn)}^{\frac{2}{n-2}}
+\|\tilde{w}\|_{L^{\critbordo}(\d\Rn)}^{\frac{2}{n-2}}\right\}\|w-\tilde{w}\|_{L^{\critbordo}(\d\Rn)}\notag
\end{align}
for all $w,\tilde{w}\in\esp_{\pares}$. 
Hence, it follows from the estimate (\ref{estim:K}) that
$$
\|\Phi_{\pares}(w)-\Phi_{\pares}(\tilde{w})\|_{\esp}
\leq C\left(\|w\|_{\esp}^{\frac{2}{n-2}}
+\|\tilde{w}\|_{\esp}^{\frac{2}{n-2}}\right)\|w-\tilde{w}\|_{\esp} 
$$
for any $w,\tilde{w}\in\esp_{\pares}$. Thus, for $\a_1$ small, the contraction maximum principle implies that the mapping $\Phi_{\pares}$ has a fixed point $\w$. Now the result follows from choosing $\v=\u+\w$. Observe that $\v$ cannot be identically zero because of (\ref{Propo5:1}) and Proposition \ref{Propo1} with $\a=\a_1$ small.
\ep

Given a pair $\pares\in\esppares$ we define
\ba\label{eq:def:energia}
\mathcal{F}_g(\xiup,\epsilon)=&\int_{\Rn}(|d\v|_g^2+c_nR_g\v^2)+\int_{\d\Rn}d_n\cmedia_g\v^2
\\
&-\frac{(n-2)^2}{n-1}\int_{\d\Rn}|\v|^{\frac{2(n-1)}{n-2}}
-\frac{n-2}{n-1}\int_{\d\Rn}\u^{\frac{2(n-1)}{n-2}}\,.\notag
\end{align}
\begin{proposition}\label{Propo6}
Suppose that $|h(x)|+|\d h(x)|+|\d^2h(x)|\leq \a_1$ for all $x\in\Rn$, where $\a_1$ is the constant obtained in Proposition \ref{Propo5}. Choosing $\a_1$ smaller if necessary, the function $\mathcal{F}_g$ is continuously differentiable and,
if $(\bar{\xiup},\bar{\epsilon})$ is a critical point of $\mathcal{F}_g$, 
then $v_{(\bar{\xiup},\bar{\epsilon})}$ is a positive smooth solution of
\begin{equation}\label{Propo6:1}
\begin{cases}
\Delta_g v_{(\bar{\xiup},\bar{\epsilon})}-c_nR_gv_{(\bar{\xiup},\bar{\epsilon})}=0\,,\:\:\:&\text{in}\:\Rn\,,
\\
\frac{\d}{\d x_n}v_{(\bar{\xiup},\bar{\epsilon})}-d_n\cmedia_g v_{(\bar{\xiup},\bar{\epsilon})}
+(n-2)v_{(\bar{\xiup},\bar{\epsilon})}^{\frac{n}{n-2}}=0\,,\:\:\:&\text{on}\:\d\Rn\,.
\end{cases}
\end{equation}
\end{proposition}

In the proof of Proposition \ref{Propo6} we will use the following removable singularities theorem, which is a slight modification of Proposition 2.7 of \cite{lee-parker}:
\begin{lemma}\label{extensao:sol}
Let $(M^n,g)$ be a Riemannian manifold with boundary $\partial M$. Let $x\in \partial M$ be a boundary point and $\mathcal{U}\subset M$ an open set containing $x$. 
Let $u$ be a weak solution to
\begin{equation}\notag
\begin{cases}
\Delta_g u +\phi u= 0\,,&\text{in}\:\mathcal{U}\backslash\{x\}\,
\\
\frac{\partial u}{\partial \eta} +\psi u=0\,,&\text{on}\:\mathcal{U}\cap\partial M\backslash\{x\}\,,
\end{cases}
\end{equation}
where $\phi\in L^{\frac{n}{2}}(\mathcal{U})$ and $\psi\in L^{n-1}(\mathcal{U}\cap\partial M)$.
Suppose that $u\in L^{q}(\mathcal{U})\cap L^{p}(\mathcal{U}\cap \d M)$ for some $q>\frac{n}{n-2}$ and $p>\frac{n-1}{n-2}$. 
Then $u$ is a weak solution to
\begin{equation}\notag
\begin{cases}
\Delta_g u +\phi u= 0\,,&\text{in}\:\mathcal{U}\,,
\\
\frac{\partial u}{\partial \eta} +\psi u=0\,,&\text{on}\:\mathcal{U}\cap\partial M\,.
\end{cases}
\end{equation}
\end{lemma}
\bp[Proof of Proposition \ref{Propo6}]
Given a pair $\pares\in\esppares$, by the definition of $\v$, there exist $b_a\pares\in\R$, $a=1,...,n$, such that
\ba
&\int_{\Rn}\left(<d\v,d\psi>_g+c_nR_g\v\psi\right)
+\int_{\d\Rn}\left(d_n\cmedia_g\v\psi-(n-2)|\v|^{\frac{2}{n-2}}\v\psi\right)\notag
\\
&\hspace{1cm}=\sum_{a=1}^nb_a\pares\cdot\int_{\d\Rn}\phi_{(\xiup,\e,a)}\psi\notag
\end{align}
for any $\psi\in\esp$. Hence, derivating the expression (\ref{eq:def:energia}) and observing the identity (\ref{eq:u:Q}), we obtain
$$
\frac{\d\mathcal{F}_g}{\d\e}\pares=2\sum_{a=1}^{n}b_a\pares\cdot\int_{\d\Rn}\phi_{(\xiup,\e,b)}\frac{\d }{\d \e}\v
$$
and
$$
\frac{\d\mathcal{F}_g}{\d\xiup_k}\pares=2\sum_{a=1}^{n}b_a\pares\cdot\int_{\d\Rn}\phi_{(\xiup,\e,a)}\frac{\d}{\d \xiup_k}\v\,,
\:\:\:\:\:k=1,...,n-1\,.
$$
On the other hand,
$$
\int_{\d\Rn}\phi_{(\xiup,\e,a)}(\v-\u)=0\,,
\:\:\:\:\:a=1,...,n\,,
$$
since $\v-\u\in\esp_{\pares}$. This implies
\ba
0&=\int_{\d\Rn}\frac{\d}{\d\e}\phi_{(\xiup,\e,a)}(\v-\u)
+\int_{\d\Rn}\phi_{(\xiup,\e,a)}\frac{\d}{\d\e}(\v-\u)\notag
\\
&=\int_{\d\Rn}\frac{\d}{\d\e}\phi_{(\xiup,\e,a)}(\v-\u)
+\int_{\d\Rn}\phi_{(\xiup,\e,a)}\frac{\d}{\d\e}\v
+\b(n)\delta_{an}\e^{-1}\notag
\end{align}
and
$$
0=\int_{\d\Rn}\frac{\d}{\d\xiup_k}\phi_{(\xiup,\e,a)}(\v-\u)
+\int_{\d\Rn}\phi_{(\xiup,\e,a)}\frac{\d}{\d\xiup_k}\v
-\b(n)\delta_{ak}\e^{-1}\,,
$$
where
$$
\b(n)=-\e\int_{\d\Rn}\phi_{(\xiup,\e,n)}\frac{\d}{\d\e}\u
=\e\int_{\d\Rn}\phi_{(\xiup,\e,k)}\frac{\d}{\d\xiup_k}\u> 0\,,
\:\:\:\:\:k=1,...,n-1\,.
$$
Thus
$$
-b_n\pares\b(n)
=\frac{\e}{2}\frac{\d\mathcal{F}_g}{\d\e}\pares
+\e\sum_{a=1}^{n}b_a\pares\cdot\int_{\d\Rn}\frac{\d}{\d\e}\phi_{(\xiup,\e,a)}(\v-\u)\,.
$$
Similarly,
$$
b_k\pares\b(n)
=\frac{\e}{2}\frac{\d\mathcal{F}_g}{\d\xiup_k}\pares
+\e\sum_{a=1}^{n}b_a\pares\cdot\int_{\d\Rn}\frac{\d}{\d\xiup_k}\phi_{(\xiup,\e,a)}(\v-\u)
$$
for $k=1,...,n-1$. Hence, if $(\bar{\xiup},\bar{\epsilon})$ is a critical point of $\mathcal{F}_g$, then there exists $C=C(n)$ such that
$$
\sum_{a=1}^n|b_a\paresbar|\leq C\|\vbar-\ubar\|_{L^{\critbordo}(\d\Rn)}\sum_{a=1}^n|b_a\paresbar|\,.
$$
By the estimate (\ref{estim:K}) and Propositions \ref{Propo1} and \ref{Propo5}, $\|v_{(\bar{\xiup},\bar{\epsilon})}-u_{(\bar{\xiup},\bar{\epsilon})}\|_{L^{\critbordo}(\d\Rn)}\leq CK^{\frac{1}{2}}\a_1$. Thus, choosing $\a_1$ small, we must have $b_a\paresbar=0$ for $a=1,...,n$. Hence,
\ba\label{Propo6:3}
&\int_{\Rn}\left(<d v_{(\bar{\xiup},\bar{\epsilon})},d\psi>_g+c_nR_g v_{(\bar{\xiup},\bar{\epsilon})}\psi\right)
\\
&\hspace{1cm}+\int_{\d\Rn}\left(d_n\cmedia_g v_{(\bar{\xiup},\bar{\epsilon})}\psi
-(n-2)|v_{(\bar{\xiup},\bar{\epsilon})}|^{\frac{2}{n-2}}v_{(\bar{\xiup},\bar{\epsilon})}\psi\right)
=0\notag
\end{align}
for any $\psi\in\esp$. 

Now we are going to show that $v_{(\bar{\xiup},\bar{\epsilon})}\geq 0$ on $\d\Rn$. To that end, we set $\psi=\min\{\vbar,0\}$ and use the equation (\ref{Propo6:3}) to conclude that
\ba\label{Propo6:2}
&\int_{\Rn\cap\{\vbar<0\}}\left(|d\vbar|_{g}^2
+c_nR_{g}\vbar^2\right)
\\
&\hspace{0.5cm}+\int_{\d\Rn\cap\{\vbar<0\}}d_n\cmedia_g\vbar^2
=(n-2)\int_{\d\Rn\cap\{\vbar<0\}}|\vbar|^{\frac{2(n-1)}{n-2}}\,.\notag
\end{align}

Using (\ref{Corol3:4}) with $w=\psi$ we see that 
\ba
\left(\int_{\d\Rn\cap\{\vbar<0\}}|\vbar|^{\frac{2(n-1)}{n-2}}\right)^{\frac{n-2}{n-1}}\notag
&\leq
2K\int_{\Rn\cap\{\vbar<0\}}\left(|d\vbar|^2_g+c_nR_g\vbar^2\right)\notag
\\
&\hspace{0.5cm}+2K\int_{\d\Rn\cap\{\vbar<0\}}d_n\cmedia_g\vbar^2\,.\notag
\end{align}
From this, together with (\ref{Propo6:2}), we deduce that $\vbar\geq 0$ almost everywhere on $\d\Rn$ or
$$
\left(\int_{\d\Rn\cap\{\vbar<0\}}|\vbar|^{\frac{2(n-1)}{n-2}}\right)^{\frac{1}{n-1}}
\geq \frac{1}{2K(n-2)}\,.
$$
On the other hand,
$$
\left(\int_{\d\Rn\cap\{\vbar<0\}}|\vbar|^{\frac{2(n-1)}{n-2}}\right)^{\frac{n-2}{2(n-1)}}
\leq 
\left(\int_{\d\Rn}|\vbar-\ubar|^{\frac{2(n-1)}{n-2}}\right)^{\frac{n-2}{2(n-1)}}
\leq CK^{\frac{1}{2}}\a_1\,.\notag
$$
Hence, choosing $\a_1$ sufficiently small we have $\vbar\geq 0$ on $\d\Rn$. In particular, the equation (\ref{Propo6:3}) can be written as 
$$
\int_{\Rn}\left(<d v_{(\bar{\xiup},\bar{\epsilon})},d\psi>_g+c_nR_g v_{(\bar{\xiup},\bar{\epsilon})}\psi\right)
+\int_{\d\Rn}\left(d_n\cmedia_g v_{(\bar{\xiup},\bar{\epsilon})}\psi
-(n-2)v_{(\bar{\xiup},\bar{\epsilon})}^{\frac{n}{n-2}}\psi\right)
=0\notag
$$
for any $\psi\in\esp$. By a result of Cherrier in \cite{cherrier}, $\vbar$ is smooth.

The fact that $v_{(\bar{\xiup},\bar{\epsilon})}> 0$ in $\Rn$ is just a consequence of the maximum principle, as follows. We set $\tilde{g}=\tilde{u}^{\frac{4}{n-2}}g$, where $\tilde{u}(x)=(1+|x|^2)^{\frac{2-n}{2}}$.
Observe that $\tilde{u}$ satisfies $\Delta\tilde{u}+n(n-2)\tilde{u}^{\frac{n+2}{n-2}}=0$ in $\Rn$ and 
we have
\ba
c_nR_{\tilde{g}}
&=-\tilde{u}^{-\frac{n+2}{n-2}}\Delta\tilde{u}
-\tilde{u}^{-\frac{n+2}{n-2}}(\Delta_g\tilde{u}-\Delta\tilde{u}-c_nR_g\tilde{u})\notag
\\
&\geq n(n-2)-C\tilde{u}^{-\frac{n+2}{n-2}}\left\{|h||\d^2\tilde{u}|
+|\d h||\d\tilde{u}|+(|\d^2h|+|\d h|^2)|\tilde{u}|\right\}\,.\notag
\end{align}
Using the facts that $h(x)=0$ for $|x|\geq 1$ and $|h|+|\d h|+|\d^2h|\leq C\a_1$ we can assume that $R_{\tilde{g}}>0$, by choosing $\a_1$ small . 

Let $S_+^n$ be a hemisphere of $S_{1/2}^n$. We will use the well known conformal equivalence between $S_+^n\backslash\{x_0\}$ and $\Rn$ realized by the stereographic projection, where $x_0\in\d S_+^n$. Under this equivalence, the standard metric on $S_+^n$ is written on $\Rn$ as $\tilde{u}^{\frac{4}{n-2}}\delta$, where $\delta$ is the Euclidean metric on $\Rn$. We set $\tilde{v}=\tilde{u}^{-1}v_{(\bar{\xiup},\bar{\epsilon})}$. By the properties (\ref{conf:proper}) of the operators $L_g=\Delta_g-c_nR_g$ and $B_g=\frac{\d}{\d\eta}-d_n\cmedia_g$, we have 
$$
L_{\tilde{g}}(\tilde{v})=\tilde{u}^{-\frac{n+2}{n-2}}L_{g}v_{(\bar{\xiup},\bar{\epsilon})}=0\,,
\:\:\:\:\text{in}\:\:S^n_+\,,
$$
and
$$
B_{\tilde{g}}(\tilde{v})+(n-2)\tilde{v}^{\frac{n}{n-2}}
=\tilde{u}^{-\frac{n}{n-2}}B_g\vbar+(n-2)(\tilde{u}^{-1}\vbar)^{\frac{n}{n-2}}
=0\,,
\:\:\:\:\text{on}\:\:\d S^n_+\,.
$$
To establish the last two equations, we also used Lemma \ref{extensao:sol}.

Since $R_{\tilde{g}}>0$, it follows from the maximum principle in $S_+^n$ and the Hopf Lemma that 
if $\tilde{v}\geq 0$ on $\d S_+^n$ then we have either $\tilde{v}>0$ or $\tilde{v}\equiv 0$ in $S_+^n$. The latter contradicts the last assertion of Proposition \ref{Propo5}.
Hence, $\tilde{v}\geq 0$ on $\d S_+^n$ implies that $\tilde{v}> 0$ in $S_+^n$.
Since we have proved that $\vbar\geq 0$ on $\d\Rn$, we conclude that $\vbar> 0$ in $\Rn$.
\ep

\section{An estimate for the energy of a bubble}\label{sec:estim:energy}
In this section we will show that the energy function $\mathcal{F}_g$ can be approximated by a certain auxiliary function.

We fix a multi-linear form $W:\R^n\times\R^n\times\R^n\times\R^n\to \R$ satisfying the algebraic properties of the Weyl tensor. We set 
$$
|W|^2=\sum_{a,b,c,d=1}^n(W_{acbd}+W_{adbc})^2
$$
and assume that $|W|^2>0$.
Recall that throughout this article we work with indices $1\leq i,j,k,l\leq n-1$ and $1\leq a,b,c,d\leq n$ and set $\bar{x}=(x_1,...,x_{n-1},0)\in\d\Rn$ whenever $x=(x_1,...,x_{n-1},x_{n})\in\Rn$ . 
For $x\in\Rn$ we set
$$
H_{ij}(x)=H_{ij}(\bar{x})=W_{ikjl}x^kx^l\:\:\:\text{and}\:\:\:H_{nb}(x)=0
$$
and define
$\bar{H}_{ab}(x)=f(|\bar{x}|^2)H_{ab}(x)$,
where
\begin{equation}\label{eq:f}
f(s)=\sum_{j=0}^{\g}a_js^j\,.
\end{equation}
The integer $0<d<\frac{n-6}{4}$ and the coefficients $a_0,...,a_{\g}\in\R$ will be chosen later. Observe that $H$ is symmetric, trace-free, independent of the coordinate $x_n$ and satisfies
$$
x^aH_{ab}(x)=x^iH_{ib}(x)=0=\d_aH_{ab}(x)=\d_iH_{ib}(x)\,,\:\:\:\:\text{for any}\:x\in\Rn\,.
$$

We define a Riemannian metric $g=\exp(h)$ on $\Rn$ where $h$ is a trace-free symmetric two tensor on $\Rn$ satisfying
\begin{equation}
\begin{cases}\notag
h_{ab}(x)=\mu\l^{2\g}f(\l^{-2}|\bar{x}|^2)H_{ab}(x)\,,\:\:\:\:&\text{for}\:|x|\leq \rho\,,
\\
h_{ab}(x)=0\,,\:\:\:\:&\text{for}\:|x|\geq 1\,.
\end{cases}
\end{equation}
Here, $\mu\leq 1$, $\l\leq\rho\leq 1$ and we suppose that $h_{nb}(x)=0$ for any $x\in\Rn$ and $\d_nh_{ab}(x)=0$ for any $x\in\d\Rn$.
We also assume that $|h|+|\d h|+|\d^2h|\leq \a_1$ where $\a_1$ is the constant obtained in Proposition \ref{Propo5}.
Observe that
$$
x^ah_{ab}(x)=x^ih_{ib}(x)=0=\d_ah_{ab}(x)=\d_ih_{ib}(x)\,,\:\:\:\:\text{for}\:|x|\leq\rho\,.
$$
and $h_{ab}(x)=O(\mu(\l+|x|)^{2\g+2})$. The second fundamental form of $\d\Rn$ satisfies 
$$
\pi_{ij}=\Gamma^n_{ij}=\frac{1}{2}(g_{in,j}+g_{jn,i}-g_{ij,n})=0\,.
$$ 
In particular, the mean curvature of $\d\Rn$ is given by $\cmedia_g=\frac{1}{n-1}g^{ij}\pi_{ij}=0$.

\bigskip
Using Proposition \ref{Propo5}, for each pair $\pares\in\esppares$ we choose $\v$ to be the unique element of $\esp$ such that $\v-\u\in\esp_{(\xiup,\e)}$ and
\begin{equation}\notag
\int_{\Rn}\left(<d\v,d\psi>_g+c_nR_g\v\psi\right)
-(n-2)\int_{\d\Rn}|\v|^{\frac{2}{n-2}}\v\psi=0
\end{equation}
for all $\psi\in\esp_{(\xiup,\e)}$.

\bigskip
Finally, we define $\Omega=\{\pares\in\esppares\,;\:|\xiup|<1,\,\frac{1}{2}<\e<2\}$\,.
\begin{proposition}\label{Propo5'}
For any pair $\pares\in\l\Omega$ we have the estimates
$$
\left\|\Delta_g\u-c_nR_g\u\right\|_{L^{\frac{2n}{n+2}}(\Rn)}
\leq C\mu\l^{2\g+2}+C\left(\frac{\l}{\rho}\right)^{\frac{n-2}{2}}
$$
and
$$
\left\|\Delta_g\u-c_nR_g\u+\mu\l^{2\g}f(\l^{-2}|\bar{x}|^2)H_{ij}\d_i\d_j\u\right\|_{L^{\frac{2n}{n+2}}(\Rn)}
\leq C\mu^2\l^{4\g+4}+C\left(\frac{\l}{\rho}\right)^{\frac{n-2}{2}}\notag\,.
$$
\end{proposition}
\bp
We just observe that
$$
|\Delta_g\u(x)-c_nR_g(x)\u(x)|\leq C\mu\l^{\frac{n-2}{2}}(\l+|x|)^{2\g+2-n}
$$
and
$$
|\Delta_g\u(x)-c_nR_g(x)\u(x)+\mu\l^{2\g}f(\l^{-2}|\bar{x}|^2)H_{ij}(x)\d_i\d_j\u(x)|\leq C\mu^2\l^{\frac{n-2}{2}}(\l+|x|)^{4\g+4-n}
$$
for $|x|\leq \rho$. In the last inequality we used the fact that, since $\d_ah_{ab}(x)=0$ for $|x|\leq \rho$, Lemma \ref{curv_esc} implies that $|R_g(x)|\leq |\d h(x)|^2+|h(x)||\d^2h(x)|$ for $|x|\leq \rho$.
\ep
\begin{corollary}\label{Corol6}
For any pair $\pares\in\l\Omega$ we have the estimate
$$
\|\v-\u\|_{L^{\frac{2n}{n-2}}(\Rn)}
+\|\v-\u\|_{L^{\frac{2(n-1)}{n-2}}(\d\Rn)}
\leq C\mu\l^{2\g+2}+
C\left(\frac{\l}{\rho}\right)^{\frac{n-2}{2}}\,.
$$
\end{corollary}
\bp
It follows from Proposition \ref{Propo5} and the estimate (\ref{estim:K}) that
\ba
\|\v-\u\|_{L^{\frac{2n}{n-2}}(\Rn)}+\|\v-\u\|_{L^{\frac{2(n-1)}{n-2}}(\d\Rn)}
&\leq C\left\|\Delta_g\u-c_nR_g\u\right\|_{L^{\frac{2n}{n+2}}(\Rn)}\notag
\\
&\leq C\mu\l^{2\g+2}+C\left(\frac{\l}{\rho}\right)^{\frac{n-2}{2}}\,,\notag
\end{align}
where we used Proposition \ref{Propo5'} in the last inequality.
\ep

In order to refine the estimate of Corollary \ref{Corol6}, using Proposition \ref{Propo4} with $h_{ab}=0$ we choose the function $\w$ to be the unique element of $\esp_{\pares}$ satisfying
\begin{equation}\label{def:w}
\int_{\Rn}<d\w,d\psi>-\int_{\d\Rn}n\u^{\frac{2}{n-2}}\w\psi
=-\int_{\Rn}\mu\l^{2\g}f(\l^{-2}|\bar{x}|^2)H_{ij}\d_i\d_j\u\psi
\end{equation}
for all $\psi \in \esp_{\pares}$. Observe that, since $x^iH_{ij}(x)=0$ for any $x\in\Rn$, we have $w_{(0,\e)}=0$. 
\begin{proposition}\label{Propo7}
The function $\w$ is smooth and satisfies, for any pair $\pares\in\l\Omega$,
$$
|\d^k\w(x)|\leq C\l^{\frac{n-2}{2}}\mu(\l+|x|)^{2\g+4-k-n}\,,\:\:\:\:\text{for all}\: x\in\Rn\,,\: k=0,1,2.
$$
\end{proposition}
\bp
First observe that there exist real numbers $b_a\pares$, $1\leq a\leq n$, such that $\w$ satisfies
\ba\label{Propo7:0}
&\int_{\Rn}<d\w,d\psi>-\int_{\d\Rn}n\u^{\frac{2}{n-2}}\w\psi
\\
&\hspace{1cm}=-\int_{\Rn}\mu\l^{2\g}f(\l^{-2}|\bar{x}|^2)H_{ij}\d_i\d_j\u\psi
+\sum_{a=1}^{n}b_a\pares\int_{\d\Rn}\phi_{(a,\xiup,\e)}\psi\notag
\end{align}
for all $\psi\in\esp$. Hence, it follows from standard elliptic theory that $\w$ is smooth. 

Now we are going to prove the pointwise estimates. Observe that
\begin{equation}\label{Propo7:1}
\left\|\mu\l^{2\g}f(\l^{-2}|\bar{x}|^2)H_{ij}(x)\d_i\d_j\u(x)\right\|_{L^{\conj}(\Rn)}
\leq C\mu\l^{2\g+2}\,.
\end{equation}
Then we apply Proposition \ref{Propo4} with $h_{ab}=0$ and use the estimates (\ref{estim:K}) and (\ref{Propo7:1}) to  conclude that
$$
\|\w\|_{L^{\crit}(\Rn)}+\|\w\|_{L^{\critbordo}(\d\Rn)}
\leq K^{\frac{1}{2}}\|\w\|_{\esp}
\leq C\mu\l^{2\g+2}\,.
$$
Moreover, we can use the equation (\ref{Propo7:0}) with $\psi=\phi_{(\xiup,\e,a)}$ to conclude that 
$$
\sum_{a=0}^{n}|b_a\pares|\leq C\mu\l^{2\g+2}\,.
$$ 
Hence,
$$
|\Delta\w(x)|
=\left|\mu\l^{2\g}f(\l^{-2}|\bar{x}|^2)H_{ij}(x)\d_i\d_j\u(x)\right|
\leq \mu\l^{\frac{n-2}{2}}(\l+|x|)^{2\g+2-n}\,,
$$
for all $x\in\Rn$, and
$$
\left|\frac{\d}{\d x_n}\w(x)+n\u^{\frac{2}{n-2}}\w(x)\right|
=\left|-\sum_{a=1}^{n}b_a\pares\phi_{(a,\xiup,\e)}(x)\right|
\leq \mu\l^{\frac{n}{2}}(\l+|x|)^{2\g+2-n}
$$
for all $x\in\d\Rn$.
\\\\
{\it{Claim.}} $\sup_{x\in\Rn}(\l+|x|)^{\frac{n-2}{2}}|\w(x)|\leq C\mu\l^{2\g+2}$

\vspace{0.2cm}
We fix $x_0\in\Rn$ and set $r=\frac{1}{2}(\l+|x_0|)$. Then we see that 
$$
\u^{\frac{2}{n-2}}(x)\leq Cr^{-1}\,,\:\:\:\:\text{for all}\:x\in B^+_{r}(x_0)\,,
$$
$$
\left|\frac{\d}{\d x_n}\w(x)+n\u^{\frac{2}{n-2}}\w(x)\right|
\leq C\mu\l^{\frac{n}{2}}r^{2\g+2-n}\,,\:\:\:\:\text{for all}\:x\in B^+_{r}(x_0)\cap\d\Rn
$$
and
$$
|\Delta\w(x)|\leq C\mu\l^{\frac{n-2}{2}}r^{2\g+2-n}\,,\:\:\:\:\text{for all}\:x\in B^+_{r}(x_0)\,.
$$
It follows from standard interior estimates that 
\ba
r^{\frac{n-2}{2}}|\w(x_0)|
&\leq 
C\|\w\|_{L^{\crit}(B^+_{r}(x_0))} 
+Cr^{\frac{n+2}{2}}\|\Delta\w\|_{L^{\infty}(B^+_{r}(x_0))}\notag
\\
&\hspace{1cm}+Cr^{\frac{n}{2}}\left\|\frac{\d}{\d x_n}\w+n\u^{\frac{2}{n-2}}\w \right\|_{L^{\infty}(B^+_{r}(x_0)\cap\d\Rn)}\notag
\\
&\leq 
C\mu\l^{2\g+2}
+C\mu\l^{\frac{n-2}{2}}r^{2\g+2+\frac{2-n}{2}}
+C\mu\l^{\frac{n}{2}}r^{2\g+2-\frac{n}{2}}\notag
\\
&\leq 
C\mu\l^{2\g+2}\,,\notag
\end{align}
since we are assuming that $d<\frac{n-6}{4}$.
This proves the Claim.

Since $\sup_{x\in\Rn}|x|^{\frac{n-2}{2}}|\w(x)|<\infty$, for all $x=(x_1,...,x_{n-1},x_n)\in\Rn$ we have
\ba
\w(x)&=
-\frac{1}{(n-2)\sigma_{n-2}}\int_{\Rn}\left(|x-y|^{2-n}+|\tilde{x}-y|^{2-n}\right)
\Delta\w(y)dy\notag
\\
&\hspace{1cm}-\frac{1}{(n-2)\sigma_{n-2}}
\int_{\d\Rn}\left(|x-y|^{2-n}+|\tilde{x}-y|^{2-n}\right)\frac{\d}{\d y_n}\w(y)dy\,,\notag
\end{align}
where $\tilde{x}=(x_1,...,x_{n-1},-x_n)$.
Now we use a bootstrap argument to prove the pointwise estimates. 
It follows from the last two inequalities that
\ba
\sup_{x\in\Rn}(\l+|x|)^{\b}|\w(x)|
&\leq
C\sup_{x\in\Rn}(\l+|x|)^{\b+2}|\Delta\w(x)|\notag
\\
&\hspace{1cm}+C\sup_{x\in\d\Rn}(\l+|x|)^{\b+1}\left|\frac{\d}{\d x_n}\w(x)\right|\notag
\end{align}
for all $0<\b<n-2$. Since
$$
|\Delta\w(x)| \leq  \mu\l^{\frac{n-2}{2}}(\l+|x|)^{2\g+2-n}\,,\:\:\:\:\text{for all}\:x\in\Rn\,,
$$
and
$$
\left|\frac{\d}{\d x_n}\w(x)\right| 
\leq  n\u^{\frac{2}{n-2}}(x)|\w(x)|
+\mu\l^{\frac{n}{2}}(\l+|x|)^{2\g+2-n}\,,\:\:\:\:\text{for all}\:x\in\d\Rn\,,
$$
we see that
$$
\sup_{x\in\Rn}(\l+|x|)^{\b}|\w(x)|
\leq
C\l \sup_{x\in\d\Rn}(\l+|x|)^{\b-1}|\w(x)|
+C\mu\l^{\b+2\g+3-\frac{n}{2}}
$$
for all $0<\b\leq n-4-2\g$. Interating we obtain
$$
\sup_{x\in\Rn}(\l+|x|)^{n-2\g-4}|\w(x)|
\leq C\mu\l^{\frac{n-2}{2}}\,.
$$
The derivative estimates follow from elliptic theory, finishing the proof.
\ep
\begin{corollary}\label{Corol8}
For any $\pares\in\l\Omega$, the function $\v-\u-\w$ satisfies
\ba
&\|\v-\u-\w\|_{L^{\frac{2n}{n-2}}(\Rn)}
+\|\v-\u-\w\|_{L^{\frac{2(n-1)}{n-2}}(\d\Rn)}\notag
\\
&\hspace{1cm}\leq C\mu^{\frac{n}{n-2}}\l^{\frac{(2\g+2)\cdot n}{n-2}}+C\left(\frac{\l}{\rho}\right)^{\frac{n-2}{2}}\,.\notag
\end{align}
\end{corollary}
\bp
It follows from the definition of $\w$ that
\ba
&\int_{\Rn}\left(<d\w,d\psi>_g+c_nR_g\w\psi\right)
-\int_{\d\Rn}n\u^{\frac{2}{n-2}}\w\psi\notag
\\
&\hspace{1cm}=-\int_{\Rn}\left\{\d_j\left((g^{ij}-\delta_{ij})\d_i\w\right)\psi-c_nR_g\w\psi\right\}\notag
\\
&\hspace{1.5cm}-\int_{\Rn}\mu\l^{2\g}f(\l^{-2}|\bar{x}|^2)H_{ij}\d_i\d_j\u\psi\,,\notag
\end{align}
for any $\psi\in\esp_{\pares}$. Hence we can write $\w=-\opinv(B_1+B_2,0)$, where
\ba
B_1&=\d_j\left((g^{ij}-\delta_{ij})\d_i\w\right)-c_nR_g\w\,,\notag
\\
B_2&=\mu\l^{2\g}f(\l^{-2}|\bar{x}|^2)H_{ij}\d_i\d_j\u\notag
\end{align}
and $\opinv$ is the operator defined in the proof of Proposition \ref{Propo5}.

On the other hand,
\ba
&\int_{\Rn}\left\{<d(\v-\u),d\psi>_g+c_nR_g(\v-\u)\psi\right\}\notag
\\
&\hspace{1cm}-\int_{\d\Rn}n\u^{\frac{2}{n-2}}(\v-\u)\psi\notag
\\
&\hspace{0.5cm}=-\int_{\Rn}\left\{<d\u,d\psi>_g+c_nR_g\u\psi\right\}\notag
\\
&\hspace{1cm}+\int_{\d\Rn}\left\{(n-2)|\v|^{\frac{2}{n-2}}\v\psi-n\u^{\frac{2}{n-2}}(\v-\u)\psi\right\}\notag
\\
&\hspace{0.5cm}=\int_{\Rn}(\Delta_g\u-c_nR_g\u)\psi
\notag
\\
&\hspace{1cm}+(n-2)\int_{\d\Rn}\left\{|\v|^{\frac{2}{n-2}}\v-\u^{\frac{n}{n-2}}
-\frac{n}{n-2}\u^{\frac{2}{n-2}}(\v-\u)\right\}\psi\,.\notag
\end{align}
Hence we can write $\v-\u=\opinv(B_3,(n-2)B_4)$, where
\ba
B_3&=\Delta_g\u-c_nR_g\u\,,\notag
\\
B_4&=(|\v|^{\frac{2}{n-2}}\v-\u^{\frac{n}{n-2}})
-\frac{n}{n-2}\u^{\frac{2}{n-2}}(\v-\u)\,.\notag
\end{align}

Puting this facts together we conclude that
$$
\v-\u-\w=\opinv(B_1+B_2+B_3,(n-2)B_4)\,.
$$

Now we are going to estimate the terms $B_1, B_2, B_3, B_4$. Since 
\ba
|B_1(x)|
&\leq C\d(|h||\d\w|)(x)+C(|\d^2h||h|+|\d h|^2)|\w|(x)\notag
\\
&\leq C\mu^{2}\l^{\frac{n-2}{2}}(\l+|x|)^{4\g+4-n}\,,\:\:\:\:\:\:\text{for}\: |x|\leq \rho\,,\notag
\end{align}
we have  
$$
\|B_1\|_{L^{\conj}(\Rn)}
\leq C\mu^{2}\l^{4\g+4}+C\left(\frac{\l}{\rho}\right)^{\frac{n-2}{2}}\,.
$$
It follows from Proposition \ref{Propo5'} that
$$
\|B_2+B_3\|_{L^{\conj}(\Rn)}
\leq C\mu^{2}\l^{4\g+4}+C\left(\frac{\l}{\rho}\right)^{\frac{n-2}{2}}\,.
$$
Since $|B_4(x)|\leq C|\v(x)-\u(x)|^{\frac{n}{n-2}}$ for any $x\in\d\Rn$, we have
$$
\|B_4\|_{L^{\conjbordo}(\d\Rn)}
\leq C\|\v-\u\|_{L^{\critbordo}(\d\Rn)}^{\frac{n}{n-2}}
\leq C\mu^{\frac{n}{n-2}}\l^{\frac{(2\g+2)\,n}{n-2}}+C\left(\frac{\l}{\rho}\right)^{\frac{n}{2}}\,,
$$
where in the last inequality we used Corollary \ref{Corol6}.

Using the estimates above we see that
$$
\|B_1+B_2+B_3\|_{L^{\conj}(\Rn)}
+\|(n-2)B_4\|_{L^{\conjbordo}(\d\Rn)}
\leq C\mu^{\frac{n}{n-2}}\l^{\frac{(2\g+2)\,n}{n-2}}+C\left(\frac{\l}{\rho}\right)^{\frac{n-2}{2}}\,.
$$
Hence,
$$
\|\v-\u-\w\|_{L^{\crit}(\Rn)}
+\|\v-\u-\w\|_{L^{\critbordo}(\d\Rn)}
\leq C\mu^{\frac{n}{n-2}}\l^{\frac{(2\g+2)\,n}{n-2}}+C\left(\frac{\l}{\rho}\right)^{\frac{n-2}{2}}\,.
$$
\ep
\begin{lemma}\label{Propo10}
For any $\pares\in\l\Omega$ we have the estimate
\ba
&\left|\int_{\d\Rn}(|\v|^{\frac{2}{n-2}}-\u^{\frac{2}{n-2}})\u\v
-\frac{1}{n-1}\int_{\d\Rn}(|\v|^{\frac{2(n-1)}{n-2}}-\u^{\frac{2(n-1)}{n-2}})\right|\notag
\\
&\hspace{1cm}\leq C\mu^{\frac{2(n-1)}{n-2}}\l^{\frac{(4\g+4)(n-1)}{n-2}}
+C\left(\frac{\l}{\rho}\right)^{n-1}\,.\notag
\end{align}
\end{lemma}
\bp 
It follows from the pointwise estimate
$$
\left|\left(|\v|^{\frac{2}{n-2}}-\u^{\frac{2}{n-2}}\right)\u\v
-\frac{1}{n-1}\left(|\v|^{\frac{2(n-1)}{n-2}}-\u^{\frac{2(n-1)}{n-2}}\right)\right|
\leq C|\v-\u|^{\frac{2(n-1)}{n-2}}
$$
that
\ba
&\left|\int_{\d\Rn}(|\v|^{\frac{2}{n-2}}-\u^{\frac{2}{n-2}})\u\v
-\frac{1}{n-1}\int_{\d\Rn}(|\v|^{\frac{2(n-1)}{n-2}}-\u^{\frac{2(n-1)}{n-2}})\right|\notag
\\
&\hspace{2cm}\leq C\left\|\v-\u\right\|_{L^{\frac{2(n-1)}{n-2}}(\d\Rn)}^{\frac{2(n-1)}{n-2}}
\leq C\left(\mu\l^{2\g+2}+\left(\frac{\l}{\rho}\right)^{\frac{n-2}{2}}\right)^{\frac{2(n-1)}{n-2}},\notag
\end{align}
where in the last inequality we used Corollary \ref{Corol6}. Now the result follows.
\ep
\begin{proposition}\label{Corol12}
Let $\mathcal{F}_g$ be the function defined by the formula (\ref{eq:def:energia}). For any pair $\pares\in\l\Omega$ we have the estimate
\ba
&\Big{|}\mathcal{F}_g\pares
-\frac{1}{2}\int_{B^+_{\rho}(0)}h_{il}h_{jl}\d_i\u\d_j\u
+\frac{c_n}{4}\int_{B^+_{\rho}(0)}(\d_l h_{ij})^2\u^2\notag
\\
&\hspace{2cm}-\int_{\Rn}\mu\l^{2\g}f(\l^{-2}|\bar{x}|^2)H_{ij}\d_i\d_j\u\,\w\Big{|}\notag
\\
&\hspace{1cm}\leq C\mu^{\frac{2(n-1)}{n-2}}\l^{\frac{(4\g+4)(n-1)}{n-2}}+C\mu\l^{2\g+2}\left(\frac{\l}{\rho}\right)^{\frac{n-2}{2}}
+C\left(\frac{\l}{\rho}\right)^{n-2}\,.\notag
\end{align}
\end{proposition}
\bp
It follows from the definition of $\v$ that
\ba
\int_{\Rn}&\left\{<d\v,d(\v-\u)>_g+c_nR_g\v(\v-\u)\right\}\notag
\\
&\hspace{1cm}-(n-2)\int_{\d\Rn}|\v|^{\frac{2}{n-2}}\v(\v-\u)=0\notag
\end{align}
Thus,
\ba\label{conseq:def:v}
&\int_{\Rn}\left\{|d\v|^2_g-<d\v,d\u>_g+c_nR_g(\v^2-\u\v)\right\}
\\
&\hspace{2cm}-(n-2)\int_{\d\Rn}\left\{|\v|^{\frac{2(n-1)}{n-2}}-|\v|^{\frac{2}{n-2}}\v\u\right\}=0\,.\notag
\end{align}
We set
\ba\label{eq:def:rho}
\varrhoup
&=\int_{\Rn}\left\{<d\u,d(\v-\u)>_g+c_nR_g\u(\v-\u)\right\}
\\
&\hspace{0.5cm}-\int_{\Rn}h_{ij}\d_i\d_j \u(\v-\u)-(n-2)\int_{\d\Rn}\u^{\frac{n}{n-2}}(\v-\u)\notag
\end{align}
Thus,
\ba\label{conseq:def:u}
\varrhoup=&\int_{\Rn}\left\{-|d\u|^2_g+<d\u,d\v>_g+c_nR_g(\u\v-\u^2)\right\}\notag
\\
&\hspace{0.5cm}-\int_{\Rn}h_{ij}\d_i\d_j\u(\v-\u)
-(n-2)\int_{\d\Rn}\left\{\u^{\frac{n}{n-2}}\v-\u^{\frac{2(n-1)}{n-2}}\right\}
\notag
\end{align}
Hence, summing (\ref{conseq:def:v}) and (\ref{conseq:def:u}),
\ba
\varrhoup=&\int_{\Rn}\left\{|d\v|^2_g
-|d\u|^2_g
+c_nR_g(\v^2-\u^2)\right\}
\\
&\hspace{0.5cm}-\int_{\Rn}h_{ij}\d_i\d_j\u(\v-\u)\notag
\\
&\hspace{0.5cm}-(n-2)\int_{\d\Rn}\left\{(|\v|^{\frac{2(n-1)}{n-2}}-\u^{\frac{2(n-1)}{n-2}})
+(\u^{\frac{2}{n-2}}-|\v|^{\frac{2}{n-2}})\,\u\v\right\}\,.\notag
\end{align}
Then
\ba
\varrhoup&=
\int_{\Rn}\left\{|d\v|^2_g+c_nR_g\v^2\right\}
-\int_{\d\Rn}\left\{\frac{(n-2)^2}{n-1}|\v|^{\frac{2(n-1)}{n-2}}
+\frac{n-2}{n-1}\u^{\frac{2(n-1)}{n-2}}\right\}
\\
&\hspace{0.5cm}-\int_{\d\Rn}\left\{\frac{n-2}{n-1}|\v|^{\frac{2(n-1)}{n-2}}
+\frac{(n-2)^2}{n-1}\u^{\frac{2(n-1)}{n-2}}\right\}\notag
\\
&\hspace{0.5cm}-(n-2)\int_{\d\Rn}(\u^{\frac{2}{n-2}}-|\v|^{\frac{2}{n-2}})\,\u\v
+2(n-2)\int_{\d\Rn}\u^{\frac{2(n-1)}{n-2}}\notag
\\
&\hspace{0.5cm}
-\int_{\Rn}\left\{|d\u|^2_g+c_nR_g\u^2
+h_{ij}\d_i\d_j\u(\v-\u)\right\}\,.\notag
\end{align}
We set
\ba
B&=\int_{\Rn}\left\{|d\u|^2_g-|d\u|^2+c_nR_g\u^2+h_{ij}\d_i\d_j\u(\v-\u)\right\}\notag
\end{align}
and observe that $\int_{\Rn}|d\u|^2=(n-2)\int_{\d\Rn}\u^{\frac{2(n-1)}{n-2}}$. Hence,
\ba\label{dif:calF:B}
\mathcal{F}_g\pares-B
&=\frac{n-2}{n-1}\int_{\d\Rn}\left\{|\v|^{\frac{2(n-1)}{n-2}}-\u^{\frac{2(n-1)}{n-2}}\right\}
\\
&\hspace{1cm}-(n-2)\int_{\d\Rn}(|\v|^{\frac{2}{n-2}}-\u^{\frac{2}{n-2}})\,\u\v
+\varrhoup\notag
\\
&= O\left\{\l^{\frac{(4\g+4)(n-1)}{n-2}}\mu^{\frac{2(n-1)}{n-2}}+\left(\frac{\l}{\rho}\right)^{n-1}\right\}
+\varrhoup\notag
\end{align}
where in the last inequality we used Lemma \ref{Propo10}.

On the other hand,
\ba\label{rel:F:B}
B&=\frac{1}{2}\int_{B^+_{\rho}(0)}h_{il}h_{jl}\d_i\u\d_j\u
-\frac{c_n}{4}\int_{B^+_{\rho}(0)}(\d_l h_{ij})^2\u^2
\\
&+\int_{\Rn}\mu\l^{2\g}f(\l^{-2}|\bar{x}|^2)H_{ij}\d_i\d_j\u\,\w+e_1+e_2+e_3+e_4+e_5\notag
\end{align}
where
\ba
e_1
&=-\int_{\Rn}h_{ij}\d_i\u\d_j\u+c_n\int_{\Rn}\d_i\d_jh_{ij}\,\u^2\,,\notag
\\
e_2
&=\int_{\Rn}(g^{ij}-\delta_{ij}+h_{ij})\d_i\u\d_j\u
-\int_{B^+_{\rho}(0)}\frac{1}{2}h_{il}h_{jl}\d_i\u\d_j\u\notag
\\
&=\int_{\Rn\backslash B^+_{\rho}(0)}(g^{ij}-\delta_{ij}+h_{ij})\d_i\u\d_j\u\notag
\\
&\hspace{1cm}+\int_{B^+_{\rho}(0)}\left\{g^{ij}-\delta_{ij}+h_{ij}-\frac{1}{2}h_{il}h_{jl}\right\}\d_i\u\d_j\u\,,\notag
\\
e_3
&=c_n\int_{\Rn}(R_g-\d_i\d_jh_{ij})\,\u^2+c_n\int_{B^+_{\rho}(0)}\frac{1}{4}(\d_lh_{ij})^2\u^2\notag
\\
&=c_n\int_{\Rn\backslash B^+_{\rho}(0)}(R_g-\d_i\d_jh_{ij})\,\u^2
+c_n\int_{B^+_{\rho}(0)}\left\{R_g+\frac{1}{4}(\d_lh_{ij})^2\right\}\u^2\,,\notag
\\
e_4&=\int_{\Rn}h_{ij}\d_i\d_j\u(\v-\u-\w)\,,\notag
\\
e_5&=-\int_{\Rn\backslash B^+_{\rho}(0)}\mu\l^{2\g}f(\l^{-2}|\bar{x}|^2)H_{ij}\d_i\d_j\u\w
+\int_{\Rn\backslash B^+_{\rho}(0)}h_{ij}\d_i\d_j\u\,\w\,.\notag
\end{align}
For the expression of $e_3$ we used the fact that $\d_jh_{ij}(x)=0$ for $|x|\leq \rho$. We are going to use this same fact in the rest of this proof. 
 
Now we are going to estimate the terms $e_1,...,e_5$. 
First observe that for $|x|\leq \rho$ we have
\begin{align}\label{estim:gij}
\left|g^{ij}(x)-\delta_{ij}+h_{ij}(x)-\frac{1}{2}h_{ij}h_{jl}(x)\right|
\leq C|h(x)|^3
&\leq C\mu^3(\l+|x|)^{6\g+6}
\\
&\leq C\mu^3(\l+|x|)^{\frac{n-1}{n-2}(4\g+4)}\notag
\end{align}
and
\begin{align}\label{estim:R}
\left|R_g(x)+\frac{1}{4}(\d_lh_{ij})^2(x)\right|
&\leq C|h(x)|^2|\d^2h(x)|+C|h(x)||\d h(x)|^2
\\
&\leq C\mu^3(\l+|x|)^{6\g+4}
\leq C\mu^3(\l+|x|)^{\frac{n-1}{n-2}(4\g+4)-2}\,.\notag
\end{align}
Here, we used Lemma \ref{curv_esc}.

From the identity $\u\d_i\d_j\u-\frac{n}{n-2}\d_i\u\d_j\u=-\frac{1}{n-2}|d\u|^2\delta_{ij}$ and the fact that $\sum_{j=1}^{n-1}h_{jj}=0$ we see that
\ba
\frac{n}{n-2}\int_{\Rn}h_{ij}\d_i\u\d_j\u
&=\int_{\Rn}h_{ij}\u\d_i\d_j\u\notag
\\
&=-\int_{\Rn}\d_ih_{ij}\,\u\d_j\u-\int_{\Rn}h_{ij}\d_i\u\d_j\u\,,\notag
\end{align}
where in the last equality we integrated by parts.
Thus,
\begin{equation}\label{int:e0}
e_1=\frac{n-2}{2(n-1)}\int_{\Rn\backslash B_{\rho}^+(0)}\d_ih_{ij}\,\u\d_j\u
+c_n\int_{\Rn\backslash B_{\rho}^+(0)}\d_i\d_jh_{ij}\,\u^2\,.
\end{equation}
Then we use the identities (\ref{int:e0}), (\ref{estim:gij}) and (\ref{estim:R}) to estimate $e_1$, $e_2$ and $e_3$ respectively and conclude that
\ba\label{estim:es}
|e_1|&\leq C\rho\left(\frac{\l}{\rho}\right)^{n-2}\,,
\\
|e_2|&\leq C\left(\frac{\l}{\rho}\right)^{n-2}+C\mu^3\l^{\frac{n-1}{n-2}(4\g+4)}\,,\notag
\\
|e_3|&\leq C\rho^2\left(\frac{\l}{\rho}\right)^{n-2}+C\mu^3\l^{\frac{n-1}{n-2}(4\g+4)}\notag
\\
|e_4|&\leq C\int_{\Rn}|h||\d^2\u||\v-\u-\w|\notag
\\
&\leq C\|h\d^2\u\|_{L^{\frac{2n}{n+2}}(\Rn)}\|\v-\u-\w\|_{L^{\frac{2n}{n-2}}(\Rn)}\notag
\\
&\leq C\left\{\mu\l^{2\g+2}+\left(\frac{\l}{\rho}\right)^{\frac{n-2}{2}}\right\}\cdot\left\{\mu^{\frac{n}{n-2}}\l^{\frac{(2\g+2)\,n}{n-2}}
+\left(\frac{\l}{\rho}\right)^{\frac{n-2}{2}}\right\}\notag
\\
&\leq C\mu^{\frac{2(n-1)}{n-2}}\l^{\frac{(4\g+4)(n-1)}{n-2}}
+\mu\l^{2\g+2}\left(\frac{\l}{\rho}\right)^{\frac{n-2}{2}}
+\mu^{\frac{n}{n-2}}\l^{\frac{(2d+2)n}{n-2}}\left(\frac{\l}{\rho}\right)^{\frac{n-2}{2}}
+\left(\frac{\l}{\rho}\right)^{n-2}\,,\notag
\\
|e_5|&\leq \rho^{2d+2}\left(\frac{\l}{\rho}\right)^{n-2}\,.\notag
\end{align}

Now we are going to estimate $\varrhoup$ using its definition (equation \eqref{eq:def:rho}).
Integrating by parts and using the second equation of (\ref{eq:U}), we obtain
\ba\label{estim:varrhoup}
|\varrhoup|&\leq
\int_{\Rn}\Big|-\Delta_g\u(\v-\u)+c_nR_g\u(\v-\u)
\\
&\hspace{5cm}-h_{ij}\d_i\d_j\u(\v-\u)\Big|\notag
\\
&\leq \left\|\Delta_g\u-c_nR_g\u+h_{ij}\d_i\d_j\u\right\|_{L^{\frac{2n}{n+2}}(\Rn)}
\|\v-\u\|_{L^{\frac{2n}{n-2}}(\Rn)}\notag
\\
&\leq C\left\{\mu^2\l^{4\g+4}+\left(\frac{\l}{\rho}\right)^{\frac{n-2}{2}}\right\}
\cdot\left\{\mu\l^{2\g+2}+\left(\frac{\l}{\rho}\right)^{\frac{n-2}{2}}\right\}\notag
\\
&\leq C\mu^3\l^{6\g+6}+C\mu\l^{2\g+2}\left(\frac{\l}{\rho}\right)^{\frac{n-2}{2}}
+C\left(\frac{\l}{\rho}\right)^{n-2}\,.\notag
\end{align}
Here, we used Proposition \ref{Propo5'} and Corollary \ref{Corol6} in the second inequality.

The result now follows from (\ref{dif:calF:B}), (\ref{rel:F:B}),  (\ref{estim:es}) and (\ref{estim:varrhoup}).
\ep


\section{Finding a critical point of an auxiliary function}\label{sec:finding}

Let us follow the notations of the last section. We define
$$
F(\xiup,\epsilon)=\frac{1}{2}\int_{\Rn}\bar{H}_{il}\bar{H}_{jl}\d_i\u\d_j\u
-\frac{c_n}{4}\int_{\Rn}(\d_l\bar{H}_{ij})^2\u^2+\int_{\Rn}\bar{H}_{ij}\d_i\d_j\u\,\z
$$
where $\z$ is the unique element of $\esp_{\pares}$ that satisfies
\begin{equation}\label{eq:def:z}
\int_{\Rn}<d\z,d\psi>-\int_{\d\Rn}n\u^{\frac{2}{n-2}}\z\psi=-\int_{\Rn}\bar{H}_{ij}\d_i\d_j\u\,\psi
\end{equation}
for any $\psi\in\esp_{\pares}$. The function $\z$ is obtained in Proposition \ref{Propo4} with $h_{ab}=0$.

In this section we will show that the function $F\pares$ has a critical point, which is a strict local minimum. 
Recall that throughout this article we use indices $1\leq i,i,j,k,l,m,p,q,r,s\leq n-1$. 

Since $\bar{H}_{ab}(-x)=\bar{H}_{ab}(x)$ for any $x\in \Rn$, the function $F\pares$ satisfies $F\pares=F(-\xiup,\e)$ for all $\pares\in\esppares$. In particular, 
\begin{equation}\label{eq:Propo13}
\frac{\d}{\d\xiup_p}F(0,\e)=\frac{\d^2}{\d\e\d\xiup_p}F(0,\e)=0\,,
\:\:\:\:\:\text{for all}\:\:\e>0\,.
\end{equation}
\begin{proposition}\label{Propo14} We have
\ba
\int_{S_r^{n-2}}(\d_lH_{ij})^2(x)x^px^q
&=\frac{2\sigma_{n-2}r^{n+2}}{(n-1)(n+1)}(W_{ipjl}+W_{iljp})(W_{iqjl}+W_{iljq})\notag
\\
&+\frac{\sigma_{n-2}r^{n+2}}{(n-1)(n+1)}(W_{ikjl}+W_{iljk})^2\delta_{pq}\notag
\end{align}
and
\ba
\int_{S_r^{n-2}}(H_{ij})^2(x)x^px^q
&=\frac{2\sigma_{n-2}r^{n+4}}{(n-1)(n+1)(n+3)}(W_{ipjl}+W_{iljp})(W_{iqjl}+W_{iljq})\notag
\\
&+\frac{\sigma_{n-2}r^{n+4}}{2(n-1)(n+1)(n+3)}(W_{ikjl}+W_{iljk})^2\delta_{pq}\,.\notag
\end{align}
\end{proposition}
\bp
Observe that
$$
\int_{S_r^{n-2}}(\d_lH_{ij})^2(x)x^px^q
=\int_{S_r^{n-2}}(W_{iljr}+W_{irjl})(W_{iljm}+W_{imjl})x^rx^mx^px^q
$$
and
$$
\int_{S_r^{n-2}}(H_{ij})^2(x)x^px^q
=\int_{S_r^{n-2}}W_{ikjl}W_{irjm}x^kx^lx^rx^mx^px^q\,.
$$
Now we just need to apply Corollary \ref{Appendix:int} in the Appendix.
\ep
\begin{proposition}\label{Propo15} We have
\ba
\int_{S_r^{n-2}}(\d_l\bar{H}_{ij})^2(x)x^px^q
&=\frac{2\sigma_{n-2}r^{n+2}}{(n-1)(n+1)(n+3)}(W_{ipjl}+W_{iljp})(W_{iqjl}+W_{iljq})\notag
\\
&\hspace{1cm}\cdot\left\{(n+3)f(r^2)^2+8r^2f(r^2)f'(r^2)+4r^4f'(r^2)^2\right\}\notag
\\
&\hspace{0.5cm}+\frac{\sigma_{n-2}r^{n+2}}{(n-1)(n+1)(n+3)}(W_{ikjl}+W_{iljk})^2\delta_{pq}\notag
\\
&\hspace{1cm}\cdot\left\{(n+3)f(r^2)^2+4r^2f(r^2)f'(r^2)+2r^4f'(r^2)^2\right\}\,.\notag
\end{align}
\end{proposition}
\bp
Since
$$
\d_l\bar{H}_{ij}(x)=f(|\bar{x}|^2)\d_lH_{ij}(x)+2f'(|\bar{x}|^2)x^lH_{ij}(x)
$$
we obtain
\ba
(\d_l\bar{H}_{ij})^2(x)
&=f(|\bar{x}|^2)^2(\d_lH_{ij})^2(x)+4f(|\bar{x}|^2)f'(|\bar{x}|^2)x^l\d_lH_{ij}H_{ij}(x)
+4|\bar{x}|^2f'(|\bar{x}|^2)^2(H_{ij})^2\notag
\\
&=f(|\bar{x}|^2)^2(\d_lH_{ij})^2(x)+8f(|\bar{x}|^2)f'(|\bar{x}|^2)(H_{ij})^2(x)
+4|\bar{x}|^2f'(|\bar{x}|^2)^2(H_{ij})^2(x)\notag
\end{align}
Hence,
\ba
\int_{S_r^{n-2}}(\d_l\bar{H}_{ij})^2(x)x^px^q
&=f(|\bar{x}|^2)^2\int_{S_r^{n-2}}(\d_lH_{ij})^2(x)x^px^q\notag
\\
&\hspace{0.5cm}
+\left(8f(|\bar{x}|^2)f'(|\bar{x}|^2)+4r^2f'(|\bar{x}|^2)^2\right)
\int_{S_r^{n-2}}(H_{ij})^2(x)x^px^q\notag
\end{align}
and the result follows from Proposition \ref{Propo14}.
\ep
\begin{corollary}\label{Corol16}
We have
\ba
&\int_{S_r^{n-2}}(\d_l\bar{H}_{ij})^2(x)=\notag
\\
&\hspace{1cm}\frac{\sigma_{n-2}r^{n}}{(n-1)(n+1)}(W_{ikjl}+W_{iljk})^2
\left\{(n+1)f(r^2)^2+4r^2f(r^2)f'(r^2)+2r^4f'(r^2)^2\right\}\,.\notag
\end{align}
\end{corollary}
\bp
By Proposition \ref{Propo15},
\ba
&r^2\int_{S_r^{n-2}}(\d_l\bar{H}_{ij})^2(x)=\sum_{p=1}^{n-1}\int_{S_r^{n-2}}(\d_l\bar{H}_{ij})^2(x)(x^p)^2\notag
\\
&=\frac{2\sigma_{n-2}r^{n+2}}{(n-1)(n+1)(n+3)}(W_{ikjl}+W_{iljk})^2
\left\{(n+3)f(r^2)^2+8r^2f(r^2)f'(r^2)+4r^4f'(r^2)^2\right\}\notag
\\
&\hspace{1cm}+\frac{\sigma_{n-2}r^{n+2}}{(n-1)(n+1)(n+3)}(n-1)(W_{ikjl}+W_{iljk})^2\notag
\\
&\hspace{3cm}\cdot\left\{(n+3)f(r^2)^2+4r^2f(r^2)f'(r^2)+2r^4f'(r^2)^2\right\}\notag
\\
&=\frac{\sigma_{n-2}r^{n+2}}{(n-1)(n+1)}(W_{ikjl}+W_{iljk})^2
\left\{(n+1)f(r^2)^2+4r^2f(r^2)f'(r^2)+2r^4f'(r^2)^2\right\}\,.\notag
\end{align}
\ep
\begin{proposition}\label{Propo17}
We have
\ba
&F(0,\epsilon)=-\frac{c_n\cdot \sigma_{n-2}}{4(n-1)(n+1)}(W_{ikjl}+W_{iljk})^2\notag
\\
&\cdot\int_{0}^{\infty}\int_{0}^{\infty}r^n\left\{(n+1)f(r^2)^2+4r^2f(r^2)f'(r^2)+2r^4f'(r^2)^2\right\}
\e^{n-2}((\e+t)^2+r^2)^{2-n}drdt\,.\notag
\end{align}
\end{proposition}
\bp
It follows from symmetry arguments that $z_{(0,\e)}=0$ and
\ba
&\int_{S_r^{n-2}}\bar{H}_{il}\bar{H}_{jl}\d_i\uo\d_j\uo(x)\notag
\\
&\hspace{1cm}=\int_{S_r^{n-2}}\frac{(n-2)^2\e^{n-2}}{((\e+x_n)^2+|\bar{x}|^2)^n}
f(|\bar{x}|^2)^2W_{iplq}W_{jrlm}x^ix^jx^px^qx^rx^m=0\,.\notag
\end{align}
Hence, we have
\ba
F(0,\epsilon)
&=-\frac{c_n}{4}\int_{\Rn}(\d_l\bar{H}_{ij})^2(x)\,\uo^2(x)\notag
\\
&=-\frac{c_n}{4}\int_{0}^{\infty}\int_{0}^{\infty}
\int_{S_r^{n-2}}(\d_l\bar{H}_{ij})^2(x)\,\uo^2(x)\,d\sigma_r(x)\,dr\,dx_n\,.\notag
\end{align}
The result now follows from Corollary \ref{Corol16}.
\ep
We write
\begin{equation}\notag
F(0,\epsilon)=
-\b_n\cdot\sum_{q=0}^{2\g}\a_q\int_{0}^{\infty}
\int_{0}^{\infty}r^{2q+n}\e^{n-2}((\e+t)^2+r^2)^{2-n}dr\,dt\,,
\end{equation}
where
$$
\b_n=\frac{c_n\cdot \sigma_{n-2}}{4(n-1)(n+1)}(W_{ikjl}+W_{iljk})^2\,,
$$
and define the coefficients $\a_q\in\R$ by the formula
\begin{equation}\label{eq:alpha}
\sum_{q=0}^{2\g}\a_qs^q=(n+1)f(s)^2+4sf(s)f'(s)+2s^2f'(s)^2\,.
\end{equation}
Here, $\g$ is the integer in the formula (\ref{eq:f}).
Changing variables $t'=t/\e$ and $r'=r/\e$ we obtain
$$
F(0,\epsilon)=
-\b_n\cdot\sum_{q=0}^{2\g}\a_q\e^{2q+4}\int_{0}^{\infty}\int_{0}^{\infty}
\frac{r^{2q+n}}{((1+t)^2+r^2)^{n-2}}drdt
$$
and, changing variables $r'=r/(1+t)$, 
$$
F(0,\epsilon)=
-\b_n\cdot\sum_{q=0}^{2\g}\a_q\e^{2q+4}\int_{0}^{\infty}\frac{1}{(1+t)^{n-5-2q}}dt\cdot
\int_{0}^{\infty}
\frac{r^{2q+n}}{(1+r^2)^{n-2}}dr
$$
Now, we have
$$
\int_{0}^{\infty}\frac{1}{(1+t)^{n-5-2q}}dt=\frac{1}{n-6-2q}
$$
and
$$
\int_{0}^{\infty}\frac{r^{2q+n}}{(1+r^2)^{n-2}}dr
=\left\{\prod_{j=0}^{q}\frac{n-1+2j}{n-5-2j}\right\}\cdot\int_{0}^{\infty}\frac{r^{n-2}}{(1+r^2)^{n-2}}dr\,,
$$
where we used Lemma \ref{int:partes}.
Hence, we can write
\begin{equation}\label{eq:F:I}
F(0,\epsilon)=
-\b_n\cdot I(\e^2)\cdot\int_{r=0}^{\infty}\frac{r^{n-2}}{(1+r^2)^{n-2}}dr
\end{equation}
where
\begin{equation}\label{eq:I}
I(s)=\sum_{q=0}^{2\g}\frac{\a_q}{n-6-2q}\left\{\prod_{j=0}^{q}\frac{n-1+2j}{n-5-2j}\right\}\,s^{q+2}\,.
\end{equation}

We will now turn our attention to the second order derivatives of the function $F \pares$.
\begin{proposition}\label{Propo19}
We have
\ba
\frac{\d^2}{\d\xiup_p\d\xiup_q}F(0,\e)
&=(n-2)^2\int_{\Rn}\frac{\e^{n-2}}{(\Zeo)^{n}}\bar{H}_{pl}(x)\bar{H}_{ql}(x)\notag
\\
&-\frac{(n-2)^2}{4}\int_{\Rn}\frac{\e^{n-2}}{(\Zeo)^{n}}(\d_l\bar{H}_{ij}(x))^2x^px^q\notag
\\
&+\frac{(n-2)^2}{8(n-1)}\int_{\Rn}\frac{\e^{n-2}}{(\Zeo)^{n-1}}(\d_l\bar{H}_{ij}(x))^2\delta_{pq}\notag\,.
\end{align}
\end{proposition}
\bp
The proof is the same of Proposition 21 of \cite{brendle2}.
\ep
\begin{proposition}\label{Propo20}
We have
\ba\label{Propo20:0}
&\frac{\d^2}{\d\xiup_p\d\xiup_q}F(0,\e)
\\
&=-\frac{2(n-2)^2\sigma_{n-2}}{(n-1)(n+1)(n+3)}(W_{ipjl}+W_{iljp})(W_{iqjl}+W_{iljq})\notag
\\
&\hspace{1cm}\cdot\int_{0}^{\infty}\int_{0}^{\infty}\frac{\e^{n-2}}{(\Zer)^{n}}r^{n+4}
\left(2f(r^2)f'(r^2)+r^2f'(r^2)^2\right)drdt\notag
\\
&-\frac{(n-2)^2\sigma_{n-2}}{2(n-1)(n+1)(n+3)}(W_{ikjl}+W_{iljk})^2\delta_{pq}\notag
\\
&\hspace{1cm}\cdot\int_{0}^{\infty}\int_{0}^{\infty}\frac{\e^{n-2}}{(\Zer)^{n}}r^{n+4}
\left\{2f(r^2)f'(r^2)+r^2f'(r^2)^2\right\}drdt\notag
\\
&+\frac{(n-2)^2\sigma_{n-2}}{4(n-1)^2(n+1)}(W_{ikjl}+W_{iljk})^2\delta_{pq}\notag
\\
&\hspace{1cm}\cdot\int_{0}^{\infty}
\int_{0}^{\infty}\frac{\e^{n-2}}{(\Zer)^{n-1}}r^{n+4}f'(r^2)^2drdt\,.\notag
\end{align}
\end{proposition}
\bp
It follows from Corollary \ref{Appendix:int} in the Appendix that
\ba
\int_{S^{n-2}_r}\bar{H}_{pl}\bar{H}_{ql}(x)
&=\int_{S^{n-2}_r}f(r^2)^2H_{pl}H_{ql}(x)
=f(r^2)^2\int_{S^{n-2}_r}W_{ipkl}W_{jqml}x^ix^jx^kx^m\notag
\\
&=\frac{\sigma_{n-2}}{2(n-1)(n+1)}(W_{ipkl}+W_{ilkp})(W_{iqkl}+W_{ilkq})r^{n+2}f(r^2)^2\notag\,.
\end{align}
Hence,
\ba\label{Propo20:1}
\int_{\Rn}\frac{\e^{n-2}\bar{H}_{pl}\bar{H}_{ql}(x)}{(\Zeo)^n}
&=\frac{\sigma_{n-2}}{2(n-1)(n+1)}(W_{ipkl}+W_{ilkp})(W_{iqkl}+W_{ilkq})
\\
&\hspace{1cm}\cdot\int_{0}^{\infty}\int_{0}^{\infty}\frac{\e^{n-2}r^{n+2}}{(\Zer)^n}f(r^2)^2dtdr\,.\notag
\end{align}
It follows from Proposition \ref{Propo15} that
\ba\label{Propo20:2}
&\int_{\Rn}\frac{\e^{n-2}(\d_l\bar{H}_{ij})^2(x)x^px^q}{(\Zeo)^n}
\\
&\hspace{0.5cm}=\frac{2\sigma_{n-2}}{(n-1)(n+1)(n+3)}(W_{ipjl}+W_{iljp})(W_{iqjl}+W_{iljq})\notag
\\
&\hspace{1cm}\cdot\int_{0}^{\infty}\int_{0}^{\infty}\frac{\e^{n-2}r^{n+2}}{(\Zer)^n}
\left\{(n+3)f(r^2)^2+8r^2f(r^2)f'(r^2)+4r^4f'(r^2)^2\right\}dtdr\notag
\\
&\hspace{0.5cm}+\frac{\sigma_{n-2}}{(n-1)(n+1)(n+3)}(W_{ikjl}+W_{iljk})^2\delta_{pq}\notag
\\
&\hspace{1cm}\cdot\int_{0}^{\infty}\int_{0}^{\infty}\frac{\e^{n-2}r^{n+2}}{(\Zer)^n}
\left\{(n+3)f(r^2)^2+4r^2f(r^2)f'(r^2)+2r^4f'(r^2)^2\right\}dtdr\,.\notag
\end{align}
and from Corollary \ref{Corol16} that
\ba\label{Propo20:3}
&\int_{\Rn}\frac{\e^{n-2}(\d_l\bar{H}_{ij})^2(x)\delta_{pq}}{(\Zeo)^{n-1}}=
\frac{\sigma_{n-2}}{(n-1)(n+1)}(W_{ikjl}+W_{iljk})^2\delta_{pq}
\\
&\hspace{0.5cm}\cdot\int_{0}^{\infty}\int_{0}^{\infty}\frac{\e^{n-2}r^{n}}{(\Zer)^{n-1}}
\left\{(n+1)f(r^2)^2+4r^2f(r^2)f'(r^2)+2r^4f'(r^2)^2\right\}dtdr\,.\notag
\end{align}
Observe that
\ba\label{Propo20:3a}
&\frac{r^n}{(\Zer)^{n-1}}\left\{(n+1)f(r^2)^2+4r^2f(r^2)f'(r^2)\right\}
\\
&\hspace{1cm}=\frac{2(n-1)r^{n+2}f(r^2)^2}{(\Zer)^{n}}
+\frac{d}{dr}\left\{\frac{r^{n+1}f(r^2)^2}{(\Zer)^{n-1}}\right\}.\notag
\end{align}
Substituting the equation \eqref{Propo20:3a} in the equation \eqref{Propo20:3} we obtain
\ba\label{Propo20:4}
&\int_{\Rn}\frac{\e^{n-2}(\d_l\bar{H}_{ij})^2(x)\delta_{pq}}{(\Zeo)^{n-1}}
\\
&\hspace{1cm}=\frac{2\sigma_{n-2}}{n+1}(W_{ikjl}+W_{iljk})^2\delta_{pq}
\cdot\int_{0}^{\infty}\int_{0}^{\infty}\frac{\e^{n-2}r^{n+2}f(r^2)^2}{(\Zer)^{n}}dtdr\notag
\\
&\hspace{1cm}+\frac{2\sigma_{n-2}}{(n-1)(n+1)}(W_{ikjl}+W_{iljk})^2\delta_{pq}
\cdot\int_{0}^{\infty}\int_{ 0}^{\infty}\frac{\e^{n-2}r^{n+4}f'(r^2)^2}{(\Zer)^{n-1}}dtdr\,,\notag
\end{align}
since we are assuming that $n>4\g+6$.
Now, using the equations \eqref{Propo20:1}, \eqref{Propo20:2} and \eqref{Propo20:4} in Proposition \ref{Propo19}, we obtain
\ba
&\frac{\d^2}{\d\xiup_p\d\xiup_q}F(0,\e)\notag
\\
&=\frac{(n-2)^2\sigma_{n-2}}{2(n-1)(n+1)}(W_{ipjl}+W_{iljp})(W_{iqjl}+W_{iljq})\notag
\\
&\hspace{1cm}\cdot\int_{0}^{\infty}\int_{0}^{\infty}\frac{\e^{n-2}}{(\Zer)^{n}}r^{n+2}f(r^2)^2drdt\notag
\\
&-\frac{(n-2)^2\sigma_{n-2}}{2(n-1)(n+1)(n+3)}(W_{ipjl}+W_{iljp})(W_{iqjl}+W_{iljq})\notag
\\
&\hspace{1cm}\cdot\int_{0}^{\infty}\int_{0}^{\infty}\frac{\e^{n-2}}{(\Zer)^{n}}r^{n+2}
\left\{(n+3)f(r^2)^2+8r^2f(r^2)f'(r^2)+4r^4f'(r^2)^2\right\}drdt\notag
\\
&-\frac{(n-2)^2\sigma_{n-2}}{4(n-1)(n+1)(n+3)}(W_{ikjl}+W_{iljk})^2\delta_{pq}\notag
\\
&\hspace{1cm}\cdot\int_{0}^{\infty}\int_{0}^{\infty}\frac{\e^{n-2}}{(\Zer)^{n}}r^{n+2}
\left\{(n+3)f(r^2)^2+4r^2f(r^2)f'(r^2)+2r^4f'(r^2)^2\right\}drdt\notag
\\
&+\frac{(n-2)^2\sigma_{n-2}}{4(n-1)(n+1)}(W_{ikjl}+W_{iljk})^2\delta_{pq}\notag
\\
&\hspace{1cm}\cdot\int_{0}^{\infty}\int_{0}^{\infty}\frac{\e^{n-2}}{(\Zer)^{n}}r^{n+2}f(r^2)^2drdt\notag
\\
&+\frac{(n-2)^2\sigma_{n-2}}{4(n-1)^2(n+1)}(W_{ikjl}+W_{iljk})^2\delta_{pq}\notag
\\
&\hspace{1cm}\cdot\int_{0}^{\infty}\int_{0}^{\infty}\frac{\e^{n-2}}{(\Zer)^{n-1}}r^{n+4}f'(r^2)^2drdt\notag
\end{align}
and the result follows after we cancel out some terms in the above equation.
\ep

Let us define constants $\b_q$, for $q=0,...,2\g-1$, by the following expression:
$$
\sum_{q=0}^{2\g-1}\b_q\,s^q=2f(s)f'(s)+sf'(s)^2\,.
$$
\begin{proposition}\label{Propo21}
We have
\ba\label{eq:J}
&\int_{0}^{\infty}\int_{0}^{\infty}\frac{\e^{n-2}}{(\Zer)^{n}}r^{n+4}\left(2f(r^2)f'(r^2)+r^2f'(r^2)^2\right)drdt
\\
&\hspace{1cm}=J(\e^2)\cdot
\int_{0}^{\infty}\frac{r^{n+2}}{(1+r^2)^n}dr\,,\notag
\end{align}
where
\begin{equation}\notag
J(s)=\sum_{q=0}^{2\g-1}\frac{\b_q s^{q+2}}{n-6-2q}
\cdot\left\{\prod_{j=0}^{q}\frac{n+3+2j}{n-5-2j}\right\}\,.
\end{equation}
\end{proposition}
\bp
\ba
&\int_{0}^{\infty}\int_{0}^{\infty}\frac{\e^{n-2}}{(\Zer)^{n}}r^{n+4}\left(2f(r^2)f'(r^2)+r^2f'(r^2)^2\right)drdt\notag
\\
&\hspace{1cm}=\sum_{q=0}^{2\g-1}\b_q
\int_{0}^{\infty}\int_{0}^{\infty}\frac{\e^{n-2}r^{n+4+2q}}{(\Zer)^n}drdt\notag
\\
&\hspace{1cm}=\sum_{q=0}^{2\g-1}\b_q\e^{2q+4}
\int_{0}^{\infty}\int_{0}^{\infty}\frac{r^{n+4+2q}}{(\Zr)^n}drdt\notag
\\
&\hspace{1cm}=\sum_{q=0}^{2\g-1}\b_q\e^{2q+4}
\int_{0}^{\infty}\frac{1}{(1+t)^{n-5-2q}}dt\int_{0}^{\infty}\frac{r^{n+4+2q}}{(1+r^2)^n}dr\notag
\end{align}
Now we observe that
$$
\int_{0}^{\infty}\frac{1}{(1+t)^{n-5-2q}}dt=\frac{1}{n-6-2q}
$$
and apply Lemma \ref{int:partes} to see that
$$
\int_{0}^{\infty}\frac{r^{n+4+2q}}{(1+r^2)^n}dr=
\left\{\prod_{j=0}^{q}\frac{n+3+2j}{n-5-2j}\right\}\cdot\int_{0}^{\infty}\frac{r^{n+2}}{(1+r^2)^n}dr\,.
$$
\ep

\subsection{The case $n\geq 53$}\label{subsec:case53}
In this case we choose $\g=1$ in the equation \eqref{eq:f}. Then the coeficients $\a_q$ in the equation \eqref{eq:alpha} are given by 
$$
\a_0=(n+1)\,a_0^2\,,\:\:\:\:\:\a_1=2(n+3)\,a_0\,a_1\,,\:\:\:\:\:\a_2=(n+7)\,a_1^2\,.
$$
Thus, derivating $I(s)$ in the expression \eqref{eq:I} we obtain
\ba
I'(s)
&=\sum_{q=0}^{2}\frac{(q+2)\,\a_q}{n-6-2q}\left\{\prod_{j=0}^{q}\frac{n-1+2j}{n-5-2j}\right\}s^{q+1}\notag
\\
&=\frac{2\a_0(n-1)}{(n-6)(n-5)}\cdot s+\frac{3\a_1(n-1)(n+1)}{(n-8)(n-5)(n-7)}\cdot s^2
+\frac{4\a_2(n-1)(n+1)(n+3)}{(n-10)(n-5)(n-7)(n-9)}\cdot s^3\notag
\\
&=\frac{2(n+1)(n-1)}{n-5}\left\{\frac{1}{n-6}a_0^2s+\frac{3(n+3)}{(n-8)(n-7)}a_0a_1s^2
+\frac{2(n+3)(n+7)}{(n-10)(n-7)(n-9)}a_1^2s^3\right\}\,.\notag
\end{align}
Now we choose $a_1=-1$ and define the polynomial $p_n$ by
$$
p_n(a_0)=\frac{a_0^2}{n-6}-\frac{3(n+3)\,a_0}{(n-8)(n-7)}
+\frac{2(n+3)(n+7)}{(n-10)(n-7)(n-9)}\,.
$$
Hence,
$$
I'(1)=\frac{2(n+1)(n-1)}{n-5}p_n(a_0)\,.
$$
The discriminant of $p_n$ is then given by
\ba
\text{discrim}(p_n)
&=\frac{(n+3)^2}{(n-7)^2(n-8)^2}\left\{9-\frac{8(n-7)(n-8)^2(n+7)}{(n+3)(n-6)(n-9)(n-10)}\right\}\notag
\\
&=\frac{(n+3)^2}{(n-7)^2(n-8)^2}\frac{q(n)}{(n+3)(n-6)(n-9)(n-10)}\,,\notag
\end{align}
where 
$$
q(n)=9(n+3)(n-6)(n-9)(n-10)-8(n-7)(n-8)^2(n+7)\,.
$$
Observe that
$$
q'(n)=4n^3-210n^2+2082n-5624
$$
and 
$$
q''(n)=6(2n^2-70n+347)\,.
$$
Since the roots $\frac{70\pm\sqrt{2124}}{4}$ of $q''$ are less than 53, we see that $q''(n)>0$ for $n\geq 53$. Since $q(53)=105696$ and $q'(53)=110340$, we conclude that $\text{discrim}(p_n)>0$ for $n\geq 53$. Hence, if we set 
$$
a_0=\frac{(n+3)(n-6)}{2(n-7)(n-8)}\left\{3+\sqrt{9-\frac{8(n-7)(n-8)^2(n+7)}{(n+3)(n-6)(n-9)(n-10)}}\right\}\,,
$$ 
then $s=1$ is critical point of $I(s)$. According to Proposition \ref{Appendix:I} in the Appendix, $I''(1)<0$ for $n\geq 53$.
 
Now we will handle $J(s)$, as defined in Proposition \ref{Propo21}. We have
$$
J(s)=\frac{(n+3)\,\b_0\,s^2}{(n-6)(n-5)}+\frac{(n+3)(n+5)\,\b_1\,s^3}{(n-8)(n-5)(n-7)}
$$
where
$$
\b_0=2\,a_0\,a_1\,\:\:\:\:\text{and}\:\:\:\:\b_1=3\,a_1^2\,.
$$
Hence,$$
J(s)=\frac{(n+3)\,a_1}{n-5}\left\{\frac{2\,a_0\,s^2}{n-6}+\frac{3(n+5)\,a_1\,s^3}{(n-8)(n-7)}\right\}\,.
$$
If we set $a_0$ and $a_1$ as above we have
\ba
J(1)
&=\frac{n+3}{(n-8)(n-5)(n-7)}\notag
\\
&\hspace{1cm}\cdot\left\{
6-(n+3)\sqrt{9-\frac{8(n-7)(n-8)^2(n+7)}{(n+3)(n-6)(n-9)(n-10)}}
\right\}\notag
\end{align}

According to Proposition \ref{Appendix:J} in the Appendix, $J(1)<0$ for $n\geq 53$.

From the equations (\ref{eq:Propo13}), (\ref{eq:F:I}), (\ref{Propo20:0}) and (\ref{eq:J}) and the above results we can conclude the following:
\begin{proposition}\label{Propo:min:local53} Suppose that $n\geq 53$. If we set $a_1=-1$ and
$$
a_0=\frac{(n+3)(n-6)}{2(n-7)(n-8)}\left\{3+\sqrt{9-\frac{8(n-7)(n-8)^2(n+7)}{(n+3)(n-6)(n-9)(n-10)}}\right\}\,,
$$ 
then $I'(1)=0$, $I''(1)<0$ and $J(1)<0$. In particular, the function $F\pares$ has a strict local minimum at the point $(0,1)$.
\end{proposition}

\subsection{The case $25\leq n\leq 52$}\label{subsec:case25}
In this case we choose $\g=4$ in the equation \eqref{eq:f}. The coeficients $\a_q$ in the equation \eqref{eq:alpha} are then given by 
\ba
&\a_0=(n+1)\,a_0^2\,,\notag
\\
&\a_1=2(n+3)\,a_0\,a_1\,,\notag
\\
&\a_2=2(n+5)\,a_0\,a_2+(n+7)\,a_1^2\,,\notag
\\
&\a_3=2(n+11)\,a_1\,a_2+2(n+7)\,a_0\,a_3\,,\notag
\\
&\a_4=2(n+15)\,a_1\,a_3+(n+17)\,a_2^2+2(n+9)\,a_0\,a_4\,,\notag
\\
&\a_5=2(n+23)\,a_2\,a_3+2(n+19)\,a_1\,a_4\,,\notag
\\
&\a_6=(n+31)\,a_3^2+2(n+29)\,a_2\,a_4\,,\notag
\\
&\a_7=2(n+39)\,a_3\,a_4\,,\notag
\\
&\a_8=(n+49)\,a_4^2\,.\notag
\end{align}

\bigskip
Thus, derivating $I(s)$ in the expression \eqref{eq:I} we obtain
$$
I'(s)
=\sum_{q=0}^{8}\frac{(q+2)\,\a_q}{n-6-2q}
\left\{\prod_{j=0}^{q}\frac{n-1+2j}{n-5-2j}\right\}s^{q+1}\notag\,.
$$
Now we choose $a_1=-3/5$, $a_2=1/8$, $a_3=-1/125$, $a_4=10^{-4}$ and define the polynomial $r_n$ by $r_n(a_0)=I'(1)$. Hence,
\ba
r_n(a_0)
&=\frac{2(n-1)(n+1)}{(n-6)(n-5)}\cdot a_0^2
+\left\{\sum_{q=1}^{4}\gamma_q(n)\frac{q+2}{n-6-2q}\prod_{j=0}^{q}\frac{n-1+2j}{n-5-2j}\right\}\cdot a_0\notag
\\
&\hspace{1cm}+\sum_{q=2}^{8}\delta_q(n)\frac{q+2}{n-6-2q}\prod_{j=0}^{q}\frac{n-1+2j}{n-5-2j}\,,\notag
\end{align}
where
$$
\gamma_1(n)=-\frac{6}{5}(n+3)\,,\:\:\:\:
\gamma_2(n)=\frac{n+5}{4}\,,\:\:\:\:
\gamma_3(n)=-\frac{2}{125}(n+7)\,,\:\:\:\:
\gamma_4(n)=\frac{n+9}{5000}\,,
$$
$$
\delta_2(n)=\frac{9(n+7)}{25},\:\:\:\:
\delta_3(n)=-\frac{3(n+11)}{20},\:\:\:\:
\delta_4(n)=\frac{1009n+16385}{40000},\:\:\:\:
\delta_5(n)=-\frac{53n+1207}{25000},
$$
$$
\delta_6(n)=\frac{89n+2709}{10^6},\:\:\:\:
\delta_7(n)=-\frac{n+39}{625000},\:\:\:\:
\delta_8(n)=\frac{n+49}{10^8}.
$$
Direct computations show that $\text{discrim}(r_n)>0$ for $25\leq n\leq 52$.

If we choose 
$$
a_0=\frac{(n-6)(n-5)}{4(n-1)(n+1)}\cdot\left\{-\sum_{q=1}^{4}\gamma_q(n)\frac{q+2}{n-6-2q}
\prod_{j=0}^{q}\frac{n-1+2j}{n-5-2j}+\sqrt{\text{discrim}(r_n)}\right\}
$$ 
then $s=1$ is critical point of $I(s)$. For $25\leq n\leq 52$, direct computations show that $I''(1)$ is of the form $-e_1-e_2\sqrt{e_3}$, where $e_1,e_2,e_3$ are positive rational numbers.

The function $J(s)$, defined in Proposition \ref{Propo21}, is written as 
$$
J(s)=\sum_{q=0}^{7}\frac{\b_q s^{q+2}}{n-6-2q}
\cdot\left\{\prod_{j=0}^{q}\frac{n+3+2j}{n-5-2j}\right\}\,.
$$
where
$$
\b_0=2\,a_0\,a_1\,,\:\:\:
\b_1=4\,a_0\,a_2+3\,a_1^2\,,\:\:\:
\b_2=6\,a_0\,a_3+10\,a_1\,a_2\,,\:\:\:
\b_3=8\,a_0\,a_4+14\,a_1\,a_3+8\,a_2^2\,,
$$
$$
\b_4=18\,a_1\,a_4+22\,a_2\,a_3\,,\:\:\:\:
\b_5=28\,a_2\,a_4+15\,a_3^2\,,\:\:\:\:
\b_6=38\,a_3\,a_4\,,\:\:\:\:
\b_7=24\,a_4^2\,.
$$
For $25\leq n\leq 52$, direct computations show that $J(1)$ is of the form $-e_1-e_2\sqrt{e_3}$, where $e_1,e_2,e_3$ are positive rational numbers.
From the equations (\ref{eq:Propo13}), (\ref{eq:F:I}), (\ref{Propo20:0}) and (\ref{eq:J}) and the above results we can conclude the following:
\begin{proposition}\label{Propo:min:local25} Suppose that $25\leq n\leq 52$. If $a_1=-3/5$, $a_2=1/8$, $a_3=-1/125$, $a_4=10^{-4}$ and 
$$
a_0=\frac{(n-6)(n-5)}{4(n-1)(n+1)}\cdot\left\{-\sum_{q=1}^{4}\gamma_q(n)\frac{q+2}{n-6-2q}
\prod_{j=0}^{q}\frac{n-1+2j}{n-5-2j}+\sqrt{\text{discrim}(r_n)}\right\}
$$ 
then $I'(1)=0$, $I''(1)<0$ and $J(1)<0$. In particular, the function $F\pares$ has a strict local minimum at the point $(0,1)$.
\end{proposition}

\section{Proof of the main theorem}\label{sec:proof}
In this section we will make use of the two-tensor $H$, defined on $\Rn$, the polynomial $f$ and the open set $\Omega\subset\esppares$, which were defined in Section \ref{sec:estim:energy}. As in Sections \ref{subsec:case53} and \ref{subsec:case25}, we fix $\g=1$ if $n\geq 53$ and $\g=4$ if $25\leq n\leq 52$. We set $D_r(0)=\{x\in\d\Rn\,;\:|x|<r\}$.

The basic ingredient in the proof of the Main Theorem is the following result:
\begin{proposition}\label{Propo24}
Assume that $n\geq 25$. Let $g$ be a smooth Riemannan metric on $\Rn$ expressed as $g=\exp(h)$, where $h$ is a symmetric trace-free two-tensor on $\Rn$ satisfying the following properties:
\begin{equation}
\begin{cases}
h_{ab}(x)=\mu\l^{2\g}f(\l^{-2}|\bar{x}|)H_{ab}(x)\,,&\text{for}\:|x|\leq\rho\,,
\\
h_{ab}(x)=0\,,&\text{for}\:|x|\geq 1\,,
\\
h_{nb}(x)=0\,,&\text{for}\:x\in\Rn\,,
\\
\d_nh_{ab}(x)=0\,,&\text{for}\:x\in\d\Rn\,,
\end{cases}
\end{equation}
where $a,b=1,...,n$. We also assume that 
$$
|h(x)|+|\d h(x)|+|\d^2h(x)|\leq \a\leq\a_1\,,\:\:\:\:\text{for all}\:x\in\Rn\,,
$$ 
where $\a_1$ is the constant obtained in Proposition \ref{Propo5}.

If $\a$ and $\mu^{-2}\l^{n-4\g-6}\rho^{2-n}$ are sufficiently small, then there exists a positive smooth function $v$ satisfying 
\begin{equation}\label{Propo24:1}
\begin{cases}
\Delta_g v-c_nR_g v=0\,,\:\:\:&\text{in}\:\Rn\,,
\\
\frac{\d}{\d x_n}v-d_n\cmedia_g v
+(n-2)v^{\frac{n}{n-2}}=0\,,\:\:\:&\text{on}\:\d\Rn
\end{cases}
\end{equation}
and
\begin{equation}\label{Propo24:2}
\int_{\d\Rn}v^{\critbordo}<\left(\frac{Q(B^n,\d B)}{n-2}\right)^{n-1}\,.
\end{equation}
Moreover, there exists $c=c(n)>0$ such that
\begin{equation}\label{Propo24:3}
\sup_{D_{\l}(0)}v\geq c\l^{\frac{2-n}{2}}\,.
\end{equation}
\end{proposition}
\bp
It follows from the fact that 
$$
(n+1)f(s)^2+4sf(s)f'(s)+2s^2f'(s)^2=(n-1)f(s)^2+2(f(s)+sf'(s))^2
$$ 
and Propostion \ref{Propo17} 
that $F(0,1)<0$. According to Propositions \ref{Propo:min:local53} and \ref{Propo:min:local25}, we can choose the coefficients $a_0,...,a_{\g}$ in the formula (\ref{eq:f}) such that the point $(0,1)$ is a strict local minimum of $F$. Hence, we can find an open set $\Omega'\subset\Omega$ such that $(0,1)\in\Omega'$ and 
$$
F(0,1)<\inf_{\pares\in\d\Omega'}F\pares<0\,.
$$
Observe that $u_{(\l\xiup,\l\e)}(\l x)=\l^{-\frac{n-2}{2}}\u(x)$ and $w_{(\l\xiup,\l\e)}(\l x)=\mu\l^{2\g+2-\frac{n-2}{2}}\z(x)$ for all $x\in\Rn$. Here, $\w$ and $\z$ are the functions defined by the formulas (\ref{def:w}) and (\ref{eq:def:z}) respectively. Thus, it follows from Proposition \ref{Corol12} that
$$
\left|\mathcal{F}_g(\l\xiup,\l\e)-\mu^{2}\l^{4\g+4}F\pares\right|
\leq
C\mu^{\critbordo}\l^{\frac{(4\g+4)\,(n-1)}{n-2}}
+C\mu\l^{2\g+2}\left(\frac{\l}{\rho}\right)^{\frac{n-2}{2}}
+C\left(\frac{\l}{\rho}\right)^{n-2}\notag
$$
for all $\pares\in\Omega$. Hence, 
\ba
\left|\mu^{-2}\l^{-4\g-4}\mathcal{F}_g(\l\xiup,\l\e)-F\pares\right|
&\leq
C\mu^{\frac{2}{n-2}}\l^{\frac{4\g+4}{n-2}}\notag
\\
&\hspace{1cm}+C\mu^{-1}\l^{\frac{n-4\g-6}{2}}\rho^{\frac{2-n}{2}}
+C\mu^{-2}\l^{n-4\g-6}\rho^{2-n}\notag
\end{align}
for all $\pares\in\Omega$. If $\mu^{-2}\l^{n-4\g-6}\rho^{2-n}$ is sufficiently small then we have
$$
\mathcal{F}_g(0,\l)<\inf_{\pares\in\d\Omega'}\mathcal{F}_g(\l\xiup,\l\e)<0\,.
$$
Thus we conclude that there exists a point $(\bar{\xiup},\bar{\e})\in\Omega'$ such that
$$
\mathcal{F}_g(\l\bar{\xiup},\l\bar{\e})=\inf_{\pares\in\Omega'}\mathcal{F}_g(\l\xiup,\l\e)<0\,.
$$
By Proposition \ref{Propo6}, the function $v=v_{(\l\bar{\xiup},\l\bar{\e})}$ obtained in Proposition \ref{Propo5} is a positive smooth solution to the equations (\ref{Propo24:1}).  
Hence, by the definition of $\mathcal{F}_g$ (see the formula (\ref{eq:def:energia})) and the formula (\ref{eq:u:Q}), we have
$$
\frac{n-2}{n-1}\int_{\d\Rn}v^{\critbordo}=
\frac{n-2}{n-1}\left(\frac{Q(B^n,\d B)}{n-2}\right)^{n-1}+\mathcal{F}(\l\bar{\xiup},\l\bar{\e})\,.
$$
This implies the inequality (\ref{Propo24:2}). 

In order to prove the inequality (\ref{Propo24:3}), observe that 
$$
\|v-u_{(\l\bar{\xiup},\l\bar{\e})}\|_{L^{\critbordo}(D_{\l}(0))}
\leq
\|v-u_{(\l\bar{\xiup},\l\bar{\e})}\|_{L^{\critbordo}(\d\Rn)}
\leq C\a
$$ 
by Propositions \ref{Propo1} and \ref{Propo5}. Hence,  
$$
|D_{\lambda}(0)|^{\frac{n-2}{2(n-1)}}\sup_{D_{\l}(0)}v
\geq
\|v\|_{L^{\critbordo}(D_{\l}(0))}
\geq
-C\a+\|u_{(\l\bar{\xiup},\l\bar{\e})}\|_{L^{\critbordo}(D_{\l}(0))}\,.
$$
Now, the inequality (\ref{Propo24:3}) follows from choosing $\a$ sufficiently small.
\ep

Now the Main Theorem follows from the next theorem, using the conformal equivalence between $B^n\backslash\{(0,...,0,-1)\}$ and $\Rn$ (see Lemma \ref{lemma:eq:conf}), the properties (\ref{conf:proper}) and Lemma \ref{extensao:sol}.
\begin{theorem}\label{Propo25}
Assume that $n\geq 25$. Then there exists a smooth Riemannian metric $g$ on $\Rn$ with the following properties:
\\\\
(a) $g_{ab}(x)=\delta_{ab}\:\:\:\text{for}\:\:\:|x|\geq 1/2$;
\\
(b) $g$ is not conformally flat;
\\
(c) $\d\Rn$ is totally geodesic with respect to the induced metric by $g$;
\\
(d) there exists a sequence of positive smooth functions $\{v_{\nu}\}_{\nu=1}^{\infty}$ satisfying
\begin{equation}
\begin{cases}
\Delta_g v_{\nu}-c_nR_g v_{\nu}=0\,,\:\:\:&\text{in}\:\Rn\,,
\\
\frac{\d}{\d x_n}v_{\nu}-d_n\cmedia_g v_{\nu}
+(n-2)v_{\nu}^{\frac{n}{n-2}}=0\,,\:\:\:&\text{on}\:\d\Rn\,,
\end{cases}
\end{equation}
for all $\nu$,
$$
\int_{\d\Rn}v_{\nu}^{\critbordo}<\left(\frac{Q(B^n,\d B)}{n-2}\right)^{n-1}\,,
$$
for all $\nu$, and $\sup_{D_{1}(0)}v_{\nu}\to\infty$ as $\nu\to\infty$.
\end{theorem}
\bp
Let $\chi:\R\to\R$ be a smooth cutoff function such that $\chi(t)=1$ for $t\leq 1$ and $\chi(t)=0$ for $t\geq 2$. We define the trace-free symmetric two-tensor $h$ on $\Rn$ by
$$
h_{ab}(x)=\sum_{N=N_0}^{\infty}\chi(4N^2|x-x_N|)2^{-\g N}f(2^N|\bar{x}-x_N|)H_{ab}(x-x_N)
$$ 
where $x_N=(\frac{1}{N},0,...,0)\in\d\Rn$. Observe that $h$ is smooth and satisfies $h_{an}(x)=0$ for $x\in\Rn$ and $\d_nh_{ab}(x)=0$ for $x\in\d\Rn$. If $N_0$ is sufficiently large, then $h_{ab}(x)=0$ for $|x|\geq \frac{1}{2}$ and $|h(x)|+|\d h(x)|+|\d^2 h(x)|\leq \a$ for $x\in\Rn$, with $\a$ sufficiently small as in Proposition \ref{Propo24}. Then we define the metric $g(x)=\exp(h(x))$ for $x\in\Rn$ and the result follows from Proposition \ref{Propo24}.
\ep

\renewcommand{\theequation}{A-\arabic{equation}}
\setcounter{equation}{0}
\renewcommand{\thetheorem}{A-\arabic{theorem}}
\setcounter{theorem}{0}
\section*{Appendix A}

In this section we establish some useful identities used in Section \ref{sec:finding}. They are simple computations which are performed in the Appendix of \cite{brendle2}. 
\begin{lemma}\label{int:partes}
We have $\int_0^{\infty}\frac{s^{\alpha}ds}{(1+s^2)^m}=\frac{2m-\alpha-3}{\alpha +1}\int_0^{\infty}\frac{s^{\alpha+2}ds}{(1+s^2)^{m}}$, for $\alpha +3<2m$.
\end{lemma}

\begin{proposition} We have
$$
\int_{S_r^{n-2}}p_k=\frac{r^2}{k(k+n-3)}\int_{S_r^{n-2}}\Delta p_k
$$
for every homogeneous polynomial $p_k$ of degree $k$.
\end{proposition}

\begin{corollary}\label{Appendix:int}
We have
$$
\int_{S^{n-2}}x_ix_j=\frac{\sigma_{n-2}}{n-1}\,,
$$
$$
\int_{S^{n-2}}x_ix_jx_kx_l=\frac{\sigma_{n-2}}{(n-1)(n+1)}
(\delta_{ij}\delta_{kl}+\delta_{ik}\delta_{jl}+\delta_{il}\delta_{jk})
$$
and
\ba
\int_{S^{n-2}}x_ix_jx_kx_lx_px_q
=\frac{\sigma_{n-2}}{(n-1)(n+1)(n+3)}
(&\delta_{ij}\delta_{kl}\delta_{pq}+\delta_{ij}\delta_{kp}\delta_{lq}+\delta_{ij}\delta_{kq}\delta_{lp}\notag
\\
+&\delta_{ik}\delta_{jl}\delta_{pq}+\delta_{ik}\delta_{jp}\delta_{lq}+\delta_{ik}\delta_{jq}\delta_{lp}\notag
\\
+&\delta_{il}\delta_{jk}\delta_{pq}+\delta_{il}\delta_{jp}\delta_{kq}+\delta_{il}\delta_{jq}\delta_{kp}\notag
\\
+&\delta_{ip}\delta_{jk}\delta_{lq}+\delta_{ip}\delta_{jl}\delta_{kq}+\delta_{ip}\delta_{jq}\delta_{kl}\notag
\\
+&\delta_{iq}\delta_{jk}\delta_{lp}+\delta_{iq}\delta_{jl}\delta_{kp}+\delta_{iq}\delta_{jp}\delta_{kl})\,.\notag
\end{align}

\end{corollary}

\renewcommand{\theequation}{B-\arabic{equation}}
\setcounter{equation}{0}
\renewcommand{\thetheorem}{B-\arabic{theorem}}
\setcounter{theorem}{0}
\section*{Appendix B}

In this section we establish some results used in Section \ref{subsec:case53}. The notations here are the same of that section. In particular, we fix $a_1=-1$ and
$$
a_0=\frac{(n+3)(n-6)}{2(n-7)(n-8)}\left\{3+\sqrt{9-\frac{8(n-7)(n-8)^2(n+7)}{(n+3)(n-6)(n-9)(n-10)}}\right\}\,.
$$ 
\begin{proposition}\label{Appendix:I}
We have $I''(1)<0$ for $n\geq 53$.
\end{proposition}
\bp
We are going to prove that $I''(1)<0$ for $n\geq 70$. If $25\leq n\leq 69$ the result follows from direct computations. We write
$$
a_0=\frac{(n+3)(n-6)}{2(n-7)(n-8)}\left\{3+\sqrt{9-\frac{8p_A(n)}{p_B(n)}}\right\}\,,
$$ 
where $p_A(n)=(n-7)(n-8)^2(n+7)$, $p_B(n)=(n+3)(n-6)(n-9)(n-10)$ and define
$$
q_L(n)=p_A(n)-p_B(n)
\:\:\:\:\text{and}\:\:\:\:
q_U(n)=\a p_B(n)-p_A(n)\,,
$$
where $\a=\frac{31439}{28800}$.
\\\\
{\it{Claim.}} $q_L(n)>0$ for $n\geq 9$ and $q_U(n)>0$ for $n\geq 70$.

\vspace{0.2cm}
In order to prove the Claim, first observe that the forth order terms of $q_L$ cancel out and we have
$q_L(n)=6n^3-114n^2+712n-1516$. Hence, $q_L''(n)=36n-228>0$ for $n\geq 7$, $q_L(9)=32$ and $q_L'(9)=118$. Thus, $q_L(n)>0$ for $n\geq 9$.

Now we observe that 
$$
q_U(n)=\frac{2639}{28800}n^4-\frac{115429}{14400}n^3+\frac{1207877}{9600}n^2
-\frac{282161}{400}n+\frac{218809}{160}\,.
$$
Hence, $q_U'''(n)=\frac{2639}{1200}n-\frac{115439}{2400}>0$ for $n\geq 70$, $q_U(70)=\frac{287074}{15}$, $q_U'(70)=\frac{178522037}{7200}$ and $q_U''(70)=\frac{10910017}{4800}$. Thus, $q_U(n)>0$ for $n\geq 70$, proving the Claim.

We asume that $n\geq 70$. In particular, we conclude from the Claim that $\a>\frac{p_A(n)}{p_B(n)}>1$, which implies
$$
\frac{2(n+3)(n-6)}{(n-7)(n-8)}>a_0>\frac{(n+3)(n-6)}{2(n-7)(n-8)}(3+\sqrt{9-8\a})\,.
$$
Now we use this estimate in 
$$
I''(1)=\frac{2(n+1)(n-1)}{n-5}
\left\{
\frac{a_0^2}{n-6}-\frac{6(n+3)a_0}{(n-8)(n-7)}+\frac{6(n+3)(n+7)}{(n-10)(n-7)(n-9)}
\right\}
$$
to see that
\ba
\frac{(n-5)I''(1)}{2(n+1)(n-1)}
<\frac{4(n+3)^3(n-6)}{(n-7)^2(n-8)^2}
&-\frac{3(3+\sqrt{9-8\a})(n+3)^2(n-6)}{(n-7)^2(n-8)^2}\notag
\\
&+\frac{6(n+3)(n+7)}{(n-10)(n-7)(n-9)}\,.\notag
\end{align}
This can be written as 
$$
I''(1)<\frac{2(n+3)(n+1)(n-1)\,\gamma(n)}{(n-8)^2(n-10)(n-5)(n-7)^2(n-9)}\,,
$$ 
where 
$$
\gamma(n)=-(5+3\sqrt{9-8\a})(n+3)(n-6)(n-10)(n-9)+6(n+7)(n-7)(n-8)^2\,.
$$
In order to complete our proof, we will show that $\gamma(n)<0$ under our assumption on the dimension. Observe that $\gamma(n)=-\frac{11}{20}n^4+\frac{481}{10}n^3-\frac{15099}{20}n^2+\frac{21162}{5}n-8205$. Hence $\gamma'''(n)=-\frac{66}{5}n+\frac{1443}{5}<0$ for $n\geq 70$, $\gamma(70)=-118392$, $\gamma'(70)=-\frac{744953}{5}$ and $\gamma''(70)=-\frac{136479}{10}$. Now the result follows.
\ep
\begin{proposition}\label{Appendix:J}
We have $J(1)<0$ for $n\geq 53$.
\end{proposition}
\bp
Let us assume that  $n\geq 53$. We want to show that $(n+3)\sqrt{9-\frac{8p_A(n)}{p_B(n)}}-6>0$, where we are using the polynomials $p_A$ and $p_B$ as in the proof of Proposition \ref{Appendix:I}. We set again $q_U(n)=\a p_B(n)-p_A(n)$ and choose $\a=\frac{7047}{6272}$.
\\\\
{\it{Claim.}} $q_U(n)>0$.

\vspace{0.2cm}
In order to prove the Claim, first observe that
$$
q_U(n)=\frac{775}{6272}n^4-\frac{27341}{3136}n^3+\frac{814983}{6272}n^2
-\frac{551233}{784}n+\frac{2063213}{1568}\,.
$$
Hence, $q_U'''(n)=\frac{2325}{784}n-\frac{82023}{1568}>0$ for $n\geq 53$, $q_U(53)=\frac{169857}{28}$, $q_U'(53)=\frac{20672955}{1568}$ and $q_U''(53)=\frac{5182395}{3136}$. Thus, $q_U(n)>0$ for $n\geq 53$, proving the Claim.

The Claim implies that 
$
(n+3)\sqrt{9-\frac{8p_A(n)}{p_B(n)}}>(n+3)\sqrt{9-8\a}\,,
$
which reduces the problem to prove that
\begin{equation}\label{Appendix:J:1}
(n+3)\sqrt{9-8\a}-6\geq 0\,.
\end{equation}
On the other hand, the fact that $\a=\frac{1}{8}\left\{9-\frac{36}{56^2}\right\}$ implies 
$\frac{1}{8}\left\{9-\frac{36}{(n+3)^2}\right\}\geq \a$, 
which is equivalent to the inequality (\ref{Appendix:J:1}).
\ep


\bigskip\noindent
\small{INSTITUTO DE MATEM\'{A}TICA \\UNIVERSIDADE FEDERAL FLUMINENSE \\NITER\'{O}I - RJ, BRAZIL\\E-mail addresses: {\bf{almaraz@vm.uff.br}}, {\bf{almaraz@impa.br}}


\begin{thebibliography}{}


\bibitem{almaraz1}
 Almaraz, S.:
 \newblock A compactness theorem for scalar-flat metrics on manifolds with boundary.
 \newblock Calc. Var. Partial Differential Equations, in press, DOI: 10.1007/s00526-010-0365-8 

\bibitem{almaraz2}
 Almaraz, S.:
 \newblock An existence theorem of conformal scalar-flat metrics on manifolds with boundary.
 \newblock Pacific J. Math. {\bf{248}}(1) , 1-22 (2010)

\bibitem{ambrosetti-li-malchiodi}
 Ambrosetti, A., Li, Y., Malchiodi, A.:
 \newblock On the Yamabe problem and the scalar curvature problem under boundary condtions.
 \newblock Math. Ann. {\bf{322}}(4), 667-699 (2002)

\bibitem{aubin1}
 Aubin, T.:
 \newblock \'{E}quations diff\'{e}rentielles non lin\'{e}aires et probl\'{e}me de Yamabe concernant la courbure scalaire.
 \newblock J. Math. Pures Appl. {\bf{55}}, 269-296 (1976)

\bibitem{berti-malchiodi}
Berti, M., Malchiodi, A.: 
 \newblock Non-compactness and multiplicity results for the Yamabe problem on $S^n$. 
\newblock J. Funct. Anal. {\bf{180}}(1), 210-241 (2001)

\bibitem{brendle2}
Brendle, S.:
 \newblock Blow-up phenomena for the Yamabe equation.
\newblock J. Amer. Math. Soc. {\bf{21}}(4), 951-979 (2008)

\bibitem{brendle-chen}
Brendle, S., Chen, S.:
\newblock An existence theorem for the Yamabe problem on manifolds with boundary.
\newblock ArXiv:0908.4327v2

\bibitem{brendle-marques}
Brendle, S., Marques, F.:
 \newblock Blow-up phenomena for the Yamabe equation II.
\newblock J. Differential Geom. {\bf{81}}, 225-250 (2009)

\bibitem{chen}
Chen, S.,:
\newblock Conformal Deformation to Scalar Flat Metrics with Constant Mean Curvature on the Boundary in Higher Dimensions.
\newblock ArXiv:0912.1302v2

\bibitem{cherrier}
Cherrier, P.:
\newblock Probl\`{e}mes de Neumann non lin\'{e}aires sur les vari\'{e}t\'{e}s Riemannienes.
\newblock J. Funct. Anal. {\bf{57}}, 154-206 (1984)


\bibitem{djadli-malchiodi-ahmedou1}
Djadli, Z., Malchiodi, A., Ould Ahmedou, M.:
 \newblock Prescribing scalar and boundary mean curvature on the three dimensional half sphere.
\newblock J. Geom. Anal. {\bf{13}}(2), 255-289 (2003)

\bibitem{djadli-malchiodi-ahmedou2}
Djadli, Z., Malchiodi, A., Ould Ahmedou, M.:
 \newblock The prescribed boundary mean curvature poblem on  $\Bbb B\sp 4$.
\newblock J. Differential Equations  {\bf{206}}(2), 373-398 (2004) 

\bibitem{druet1}
Druet, O.:
 \newblock From one bubble to several bubbles: the low-dimensional case.
\newblock J. Differential Geom. {\bf{63}}(3), 399-473 (2003)

\bibitem{druet2}
Druet, O.:
 \newblock Compactness for Yamabe metrics in low dimensions.
\newblock Int. Math. Res. Not. {\bf{23}}, 1143-1191 (2004) 

\bibitem{escobar2}
 Escobar, J.:
 \newblock The Yamabe problem on manifolds with boundary.
 \newblock J. Differential Geom. {\bf{35}}, 21-84 (1992)


\bibitem{escobar3}
 Escobar, J.:
 \newblock Conformal deformation of a Riemannian metric to a scalar flat metric with constant mean curvature on the boundary.
 \newblock Ann. Math. {\bf{136}}, 1-50 (1992)

\bibitem{escobar4}
 Escobar, J.:
 \newblock Conformal metrics with prescribed mean curvature on the boundary.
 \newblock Calc. Var. Partial Differential Equations {\bf{4}}, 559-592 (1996)

\bibitem{ahmedou-felli1}
 Felli, V., Ould Ahmedou, M.:
 \newblock Compactness results in conformal deformations of Riemannian metrics on manifolds with boundaries.
 \newblock Math. Z. {\bf{244}}, 175-210 (2003)

\bibitem{ahmedou-felli2}
 Felli, V., Ould Ahmedou, M.:
 \newblock A geometric equation with critical nonlinearity on the boundary.
 \newblock Pacific J. Math. {\bf{218}}(1), 75-99 (2005)


\bibitem{han-li}
 Han, Z., Li, Y.:
 \newblock The Yamabe problem on manifolds with boundary: existence and compactness results.
 \newblock Duke Math. J. {\bf{99}}(3), 489-542 (1999)

\bibitem{khuri-marques-schoen}
Khuri, M., Marques, F., Schoen, R.:
 \newblock A compactness theorem for the Yamabe problem.
\newblock J. Differential Geom. {\bf{81}}(1), 143-196 (2009)

\bibitem{lee-parker}
 Lee, J., Parker, T.:
 \newblock The Yamabe problem.
 \newblock Bull. Amer. Math. Soc. {\bf{17}}, 37-91 (1987)


\bibitem{li-zhang}
 Li, Y., Zhang, M.:
 \newblock Compactness of solutions to the Yamabe problem, II.
 \newblock Calc. Var. Partial Differential Equations {\bf{24}}, 185-237 (2005)


\bibitem{li-zhang2}
 Li, Y., Zhang, M.:
 \newblock Compactness of solutions to the Yamabe problem, III.
 \newblock J. Funct. Anal. {\bf{245}}(2), 438-474 (2007)


\bibitem{li-zhu2}
 Li, Y., Zhu, M.:
 \newblock Yamabe type equations on three dimensional Riemannian manifolds.
 \newblock Commun. Contemp. Math. {\bf{1}}(1), 1-50 (1999)

\bibitem{marques}
 Marques, F.:
 \newblock A priori estimates for the Yamabe problem in the non-locally conformally flat case.
 \newblock  J. Differential Geom. {\bf{71}}, 315-346 (2005)

\bibitem{coda1}
 Marques, F.:
 \newblock Existence results for the Yamabe problem on manifolds with boundary.
 \newblock Indiana Univ. Math. J. {\bf{54}}(6), 1599-1620 (2005)

\bibitem{coda2}
 Marques, F.:
 \newblock Conformal deformation to scalar flat metrics with constant mean curvature on the boundary.
 \newblock Comm. Anal. Geom. {\bf{15}}(2), 381-405 (2007)

\bibitem{ouldahmedou}
Ould Ahmedou, M.:
\newblock On the prescribed scalar and zero mean curvature on 3-D manifolds with umbilic boundary. \newblock Advanced Nonlinear Studies {\bf{6}}, 13-46 (2006)

\bibitem{schoen1}
 Schoen, R.:
 \newblock Conformal deformation of a Riemannian metric to constant scalar curvature.
 \newblock J. Differential Geom. {\bf{20}}, 479-495 (1984)
\bibitem{schoen4}
 Schoen, R.:
 \newblock On the number of constant scalar curvature metrics in a conformal class.
 \newblock Differential Geometry: A sysposium in honor of Manfredo Do Carmo (H.B.Lawson and K.Tenenblat, eds.), Wiley, 311-320 (1991)

\bibitem{schoen-yau3}
 Schoen, R., Yau, S.-T.:
 \newblock Lectures on Differential Geometry,
 \newblock Conference Proceedings and Lecture Notes in Geometry and Topology, I. International Press, Cambridge (1994)

\bibitem{schoen-zhang}
 Schoen, R., Zhang, D.:
 \newblock Prescribed scalar curvature on the n-sphere.
 \newblock Calc. Var. Partial Differential Equations {\bf{4}}(1), 1-25 (1996)


\bibitem{trudinger}
Trudinger, N.:
\newblock Remarks concerning the conformal deformation of a Riemannian structure on compact manifolds.
\newblock Ann. Scuola Norm. Sup. Pisa Cl. Sci. (3) {\bf{22}}, 265-274 (1968)

\bibitem{yamabe}
Yamabe, H.:
\newblock On a deformation of Riemannian structures on compact manifolds.
\newblock Osaka Math. J. {\bf{12}}, 21-37 (1960)

\end{thebibliography}
\end{document}